\documentclass{article}

\usepackage[pages=all, color=white, position={current page.south}, placement=bottom, scale=1, opacity=1, vshift=5mm]{background}

\usepackage[margin=1in]{geometry} 

\usepackage{amsmath}
\usepackage{amsthm}
\usepackage{amssymb}
\usepackage{algorithmic}
\usepackage{algorithm}
\usepackage{bbold}

\usepackage[utf8]{inputenc}

\usepackage[sort&compress,numbers,square]{natbib}
\usepackage{graphicx, color}

\usepackage{mathrsfs} 

\DeclareRobustCommand{\rchi}{{\mathpalette\irchi\relax}}
\newcommand{\irchi}[2]{\raisebox{\depth}{$#1\chi$}} 
\newcommand{\Ptau}{P_{\bar\tau} }
\newcommand{\NSS}{\mathrm{N}_\text{SS}}
\newcommand{\Ns}{\mathrm{N}_\textsc{samp}}
\newcommand{\Npc}{\mathrm{N}_\textsc{PC}}
\newcommand{\NKL}{\mathrm{N}_\textsc{KL}}

\newcommand{\V}{\mathbb{V}}
\newcommand{\E}{\mathbb{E}}

\newcommand{\ind}{\rchi}

\title{Global sensitivity analysis of rare event probabilities}
\author{Michael Merritt$^1$ \and Alen Alexanderian$^1$ \and Pierre Gremaud$^{1,2}$}

\date{
	$^1$Department of Mathematics, NC State University, Raleigh, NC 27695-8205, USA\\%
	$^2$The Graduate School and Department of Mathematics, NC State University, Raleigh, NC 27695-7102, USA\\[2ex]%
}

\begin{document}
	\maketitle
	
	\begin{abstract}
		By their very nature, rare event probabilities are expensive to
compute; they are also delicate to estimate as their value strongly depends  on
distributional assumptions on the model parameters. 
Hence, understanding the sensitivity of the computed rare event probabilities
to the hyper-parameters that define the distribution law of the model
parameters is crucial.
We show that by (i)
accelerating the calculation of rare event probabilities through subset
simulation and (ii) approximating the resulting probabilities through
a polynomial chaos expansion, the global sensitivity of such problems can be
analyzed through a double-loop sampling approach. The resulting method is
conceptually simple and computationally efficient; its performance is
illustrated on a subsurface flow application and on an analytical example.
		
	\end{abstract}

	\section{Introduction}\label{sec:Introduction}
Quantifying  rare event probabilities is often needed when  modeling under
uncertainty~\cite{ullmann2015multilevel, tong2020extreme, beckzuev, morio2011influence}. Rare
 events are commonly  associated with system failures or  anomalies which
pose a risk; it is thus imperative that rare event probabilities be  computed reliably. 
For the sake of concreteness,  we consider $q$ to be a scalar-valued quantity of interest (QoI) whose inputs are drawn
from the sample space $\Theta \subseteq \mathbb{R}^d$ 
with associated sigma algebra $\mathcal{F}$ and
probability measure $\mathbb{P}$.  For a given threshold
$\bar{\tau}$, the corresponding rare event probability is defined as
\begin{equation}\label{eq:Prob_RE} 
\Ptau = \mathbb{P}(q(\boldsymbol{\theta}) >
\bar{\tau}), 
\end{equation} 
where $\boldsymbol{\theta} \in \Theta$ 
is a random vector whose entries 
represent uncertain model parameters. 
Rare event
probabilities are notoriously challenging to compute; indeed, basic Monte-Carlo
simulations of (\ref{eq:Prob_RE}) are  inefficient in this context for the simple reason
that few samples actually hit the rare event domain. Several methods have
been proposed to compute $\Ptau$ more efficiently, ranging from importance
sampling and Taylor series approximations to subset simulation, the latter of
which we use in this article, see for instance \cite{beckzuev} and
Section~\ref{sec:rareeventsimulation}.

The evaluation of the rare event probability (\ref{eq:Prob_RE})  requires the distribution law governing the model 
parameters $\boldsymbol{\theta}$. In practice, such a law is typically \emph{assumed}.
Clearly, $\Ptau$ depends on these assumptions; should they be misguided,
the resulting rare event probability is likely be misleading. We  let $\boldsymbol{\xi}$ denote 
a set of hyper-parameters charecterizing the
distribution law of $\boldsymbol{\theta}$. It is  crucial to understand the sensitivity of 
$\Ptau$ to  $\boldsymbol\xi$. 
In this article, we develop an efficient method to quantify, through global sensitivity analysis (GSA), the robustness of $\Ptau$ to the choice of hyper-parameters characterizing the distribution law of  the model parameters.

To account for the uncertainty in  $\boldsymbol\xi$, we model the corresponding hyper-parameters as random variables. The rare event 
probability  takes the form 
\begin{equation}
\Ptau(\boldsymbol{\xi}) = 
\mathbb{P}(
\{q(\boldsymbol{\theta}) > \bar{\tau}\} \mid \boldsymbol{\xi} ). \label{eq:qoixi}
\end{equation}
A number of recent studies have considered how to  assess the sensitivity of  rare event estimation procedures to uncertain inputs and/or  to the  distributions of these inputs. There is a general consensus that the naive double-loop approach\,---\,whereby for each realization of $\boldsymbol{\xi}$ 
multiple samples of $\boldsymbol{\theta}$ are used to estimate  $\Ptau$\,---\,is infeasible but for the simplest of problems. An early work~\cite{morio2011influence} combines rare event estimation techniques with the traditional Monte Carlo approach for GSA of the hyper-parameters. Several studies introduce new sensitivity measures~\cite{chabridon2018reliability, lemaitre2015density, dupuis2020sensitivity, ehre2020framework} which are tailored to make the rare event SA process more tractable. Others perform sensitivity analysis in the joint space of both input parameters and hyper-parameters~\cite{chabridon2018reliability, ehre2020framework, wang2020augmented, wang2021global}. These methods increase computational efficiency through use of local SA methods~\cite{chabridon2018reliability}, surrogate models~\cite{ehre2020framework}, kernel density estimates~\cite{wang2020augmented}, and Kriging~\cite{wang2021global}. A thorough overview of current methods at the intersection of SA and rare event simulation can be found in~\cite{chabridonthesis} .

Our main contribution is to show that a double-loop approach can in fact be not only feasible but computationally expedient in order to perform GSA of $\Ptau(\boldsymbol{\xi})$ with respect to $\boldsymbol\xi$. This may seem counterintuitive since, while informative,  this type of  second level sensitivity analysis is expensive. Our approach is however structurally simpler than most of the previously cited work and achieves computational efficiency   through a combination of  fast methods for rare event 
simulations together with  the use of surrogate models. Specifically, 
we rely on subset simulation~\cite{beckzuev} to estimate
rare event probabilities and approximate $\Ptau(\boldsymbol{\xi})$ using
a polynomial chaos expansion (PCE), see respectively in Section~\ref{sec:rareeventsimulation} and Section~\ref{sec:SurrogatesforGSA}.  GSA is performed through a variance-based approach: crucially, the Sobol' indices for appropriate approximations to 
$\Ptau(\boldsymbol{\xi})$ can then be obtained ``for free'' through
analytical formul\ae. 
To  demonstrate the efficiency gains of the proposed 
method, we present an illustrative example in Section~\ref{sec:motivatingexample} and deploy our approach on it in Section~\ref{sec:ResultsAnalytic}. In Section~\ref{sec:SubsurfaceFlow}, we apply the method to a Darcy flow
problem requiring multiple estimates of the rare event probability to show
feasibility in a more computationally demanding framework. We discuss
additional challenges, perspectives and future work in Section~\ref{sec:Conclusion}. 

\section{A motivating example}\label{sec:motivatingexample}
We consider the following illustrative example~\cite{vsehic2020estimation,
beckzuev,papaioannou2015mcmc} throughout the article 
\begin{equation} \label{eq:limit_state_fun}
q(\boldsymbol\theta) = -\frac{1}{\sqrt{d}}\sum_{i = 1}^d \theta_i, 
\end{equation}
where $q$ is the QoI in (\ref{eq:Prob_RE}) 
and $\boldsymbol \theta = \begin{bmatrix} \theta_1 \, 
\cdots \, \theta_d\end{bmatrix}^\top$ with independent normally distributed entries $\theta_i \sim \mathcal{N}(\mu_i, \sigma_i^2)$, $i = 1, \ldots, d$. It is elementary to check that, for any values of the hyper-parameters 
$\boldsymbol\xi = \begin{bmatrix} 
\mu_1 & \dots & \mu_d & \sigma_1^2 & \dots & \sigma_d^2\end{bmatrix}^\top$
\begin{eqnarray}
q \sim \mathcal{N}(\bar\mu, \bar\sigma^2) \quad \mbox{ with } \left\{\begin{array}{l} \bar \mu = -\frac 1{\sqrt{d}} \sum_{i=1}^d \mu_i\\
\bar \sigma^2 = \frac 1d \sum_{i=1}^d \sigma_i^2.\end{array}\right. \label{eq:exactq}
\end{eqnarray}
For a given $\boldsymbol\xi$, the rare event probability is simply 
\begin{eqnarray}
\Ptau(\boldsymbol{\xi}) = \frac 12 -  \frac 12\operatorname{erf}\left( \frac{ \bar\tau - \bar\mu}{\sqrt{2} \bar\sigma} \right).
 \label{eq:exactp}
\end{eqnarray}

\begin{figure}[!ht]
\centering
\includegraphics[width=.33\textwidth]{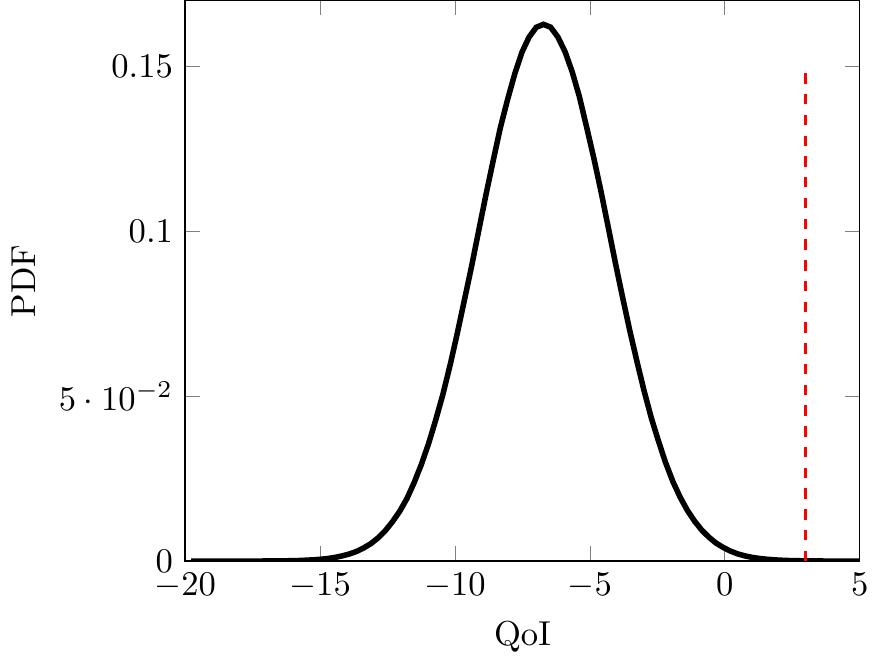}
\includegraphics[width=.33\textwidth]{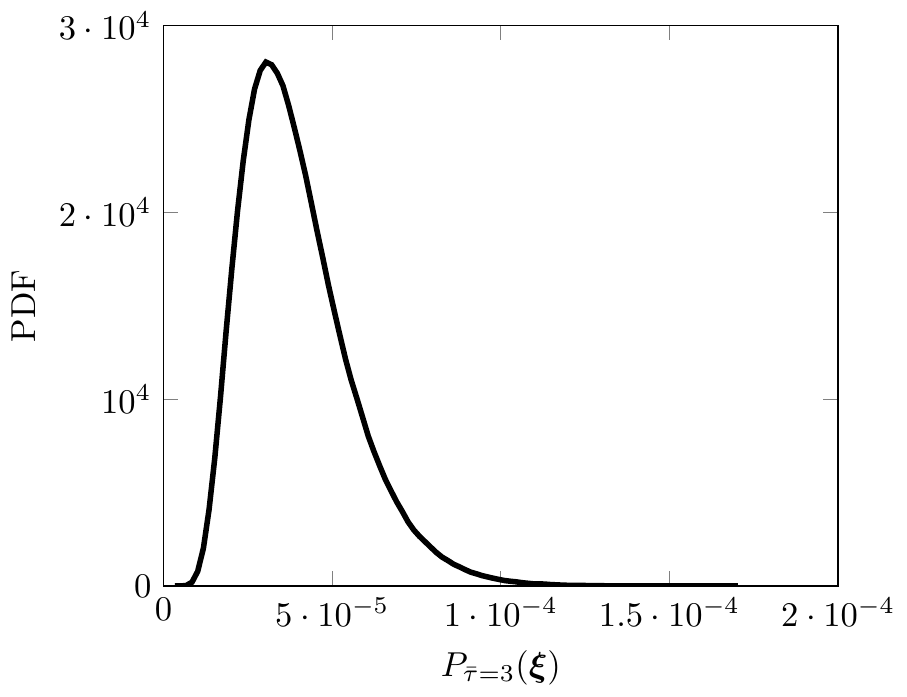}
\includegraphics[width=.33\textwidth]{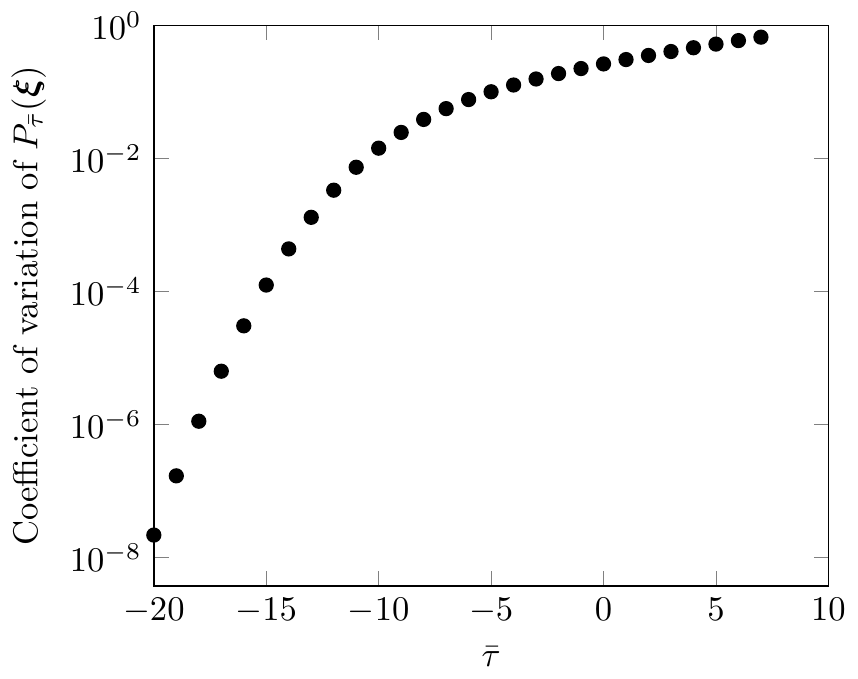}
\caption{
Left: Probability Distribution Function (PDF) of $q$ from (\ref{eq:limit_state_fun}) with the rare event threshold $\bar\tau = 3$ indicated by 
a vertical line; middle: PDF of $P_{3}(\boldsymbol\xi)$, note that from (\ref{eq:exactp}), $P_{3}(\boldsymbol\xi_{nom}) \approx 3.69\times 10^{-5}$; right: coefficient of variation of  $P_{\bar\tau}(\boldsymbol\xi)$
(ratio of standard deviation to mean) 
as $\bar{\tau}$ varies.}
\label{fig:example_problem}
\end{figure}

We model the uncertainty in the hyper-parameters by 
considering them as independent uniformly distributed random variables
with a 10 percent perturbation around their respective nominal values. 
Figure~\ref{fig:example_problem} illustrates the case  $d=5$ with $\boldsymbol\xi_{nom} = \begin{bmatrix} 1 & 2 & 3 & 4 & 5 & 10 & 8 & 6 & 4 & 2 \end{bmatrix}^\top$ as the nominal 
value for $\boldsymbol\xi$.
In particular, Figure~\ref{fig:example_problem}, right, shows how the
uncertainty in $\Ptau$ changes as $\bar\tau$ varies. As $\bar\tau$ increases,
i.e., as the event becomes rarer, the uncertainty in $\Ptau$\,---\,measured
through its coefficient of variation\,---\,increases.  We contend that this
latter behavior is generic for rare event simulations, establishing the need
for methods allowing the quantification of the effects of hyper-parameter
choices on the uncertainty in $\Ptau$.

To provide qualitative insight, we present 
a rough estimate for 
the decrease in the coefficient of variation of $\Ptau$, as the event becomes less rare.
We consider a generic $\Ptau(\boldsymbol\xi)$ as defined
in~\eqref{eq:qoixi} and assume $\Ptau$ is a random variable (i.e., 
a measurable function of $\boldsymbol\xi$). Let
$\mu = \E(\Ptau)$ and $\sigma^2 = \V(\Ptau)$ be the 
mean and variance of $\Ptau$. Recall that the coefficient
of variation of $\Ptau$ is given by $\delta(\Ptau) = \sigma/\mu$.
Note that for every $\boldsymbol\xi$, we have
$0 \leq \Ptau(\boldsymbol\xi) \leq 1$; thus,
$\Ptau(\boldsymbol\xi) \geq \Ptau(\boldsymbol\xi)^2$ and
\begin{equation}\label{equ:varbd}
    \sigma^2 = \E(\Ptau^2) - \mu^2 \leq \mu - \mu^2 = \mu(1-\mu).
\end{equation}
Therefore, 
$\delta^2(\Ptau) = \sigma^2/\mu^2 \leq \mu(1 - \mu)/\mu^2 = (1-\mu)/\mu$. 
Note that as the event becomes less rare, $\mu$ will grow 
resulting in the diminishing of the bound on the coefficient of variation. 
We point out that the inequality~\eqref{equ:varbd} can be obtained directly from 
the more general Bhatia--Davis inequality~\cite{BhatiaDavis00}.

\section{Rare event simulation}\label{sec:rareeventsimulation}
Monte Carlo simulation is a standard way of approximating the rare event probability $\Ptau$
defined in (\ref{eq:Prob_RE}). Observe that
%
%
\begin{equation}
\Ptau = \mathbb{E}[\ind_{\bar{\tau}}] = \int_{\Theta} \ind_{\bar{\tau}}
(\boldsymbol\theta) \pi(\boldsymbol{\theta}) ~d\boldsymbol{\theta}, 
\label{ptauMC}
\end{equation}
where $\ind_{\bar \tau}$ denotes the indicator function of the set 
$\{ \boldsymbol \theta \in \Theta : q(\boldsymbol\theta) > \bar\tau\}$ and 
$\pi(\boldsymbol{\theta})$ is the PDF of $\boldsymbol{\theta}$. This leads to 
the following Monte Carlo (MC) estimator 
\begin{equation}
\label{eq:Prob_RE_MC} \hat{P}_{\bar{\tau}}^{MC} = \frac{1}{N} \sum_{i=1}^N
\ind_{\bar{\tau}} (\boldsymbol{\theta}^{(i)}), 
\end{equation} 
where $\boldsymbol{\theta}^{(i)}$, $i=1,\dots,N$,  are (independent) realizations of
$\boldsymbol{\theta}$.

In the case of rare events, i.e., of small probabilities $\Ptau$, the basic MC 
estimator~(\ref{eq:Prob_RE_MC}) becomes  computationally inefficient. 
Indeed, consider the coefficient of variation $\delta
(\hat{P}_{\bar{\tau}}^{MC})$ of the above estimator and observe
\begin{equation} \label{eq:coef_var_MC} 
\delta^2\left(
\hat{P}_{\bar{\tau}}^{MC} \right) =
\frac{\V\left(\hat{P}_{\bar{\tau}}^{MC}
\right)}{\mathbb{E}\left[\hat{P}_{\bar{\tau}}^{MC} \right]^2}  = \frac{1 -
P_{\bar{\tau}}}{N P_{\bar{\tau}}} \approx \frac 1{N P_{\bar{\tau}}} \quad \mbox{if $0< P_{\bar\tau} \ll 1$}.
\end{equation} 
In other words, ensuring a given accuracy requires 
$N \approx \frac 1{\Ptau \delta^2}$. For increasingly rare events, i.e. small $\Ptau$, the error in \eqref{eq:Prob_RE_MC} will increase accordingly. Standard MC methods are thus poor candidates for rare event estimation.


\paragraph*{The subset simulation method}
We rely on the
subset simulation (SS) method~\cite{au2001estimation, schueller2004critical}
to accelerate rare event computation. 
This approach  decomposes the rare event estimation problem into a series of
``frequent event'' estimation problems that are more tractable; it has been observed that this may reduce the coefficient of variation by more than an order of magnitude over standard MC~\cite{au2001estimation, beckzuev, schueller2004critical}. This corresponds to a substantially lower computational burden for estimating
 rare event probabilities.  

Consider the \textit{rare event domain} $F = \{ \boldsymbol{\theta}
\in \Theta ~|~ q(\boldsymbol{\theta}) > \bar{\tau} \}$ and 
a sequence of nested subsets of $F$
\[
F = F_{L} \subset \dots \subset F_{2} \subset  F_{1},
\]
where $F_i = \{ \boldsymbol{\theta} \in \Theta ~|~ q(\boldsymbol{\theta}) > 
\tau_i \}$, $i = 1, \ldots, L$ with  $\tau_1 < \tau_2 < \dots < \tau_L = \bar{\tau} $. The rare event probability $\Ptau$ can thus be decomposed into a product of conditional probabilities
\begin{equation} 
\Ptau = \mathbb{P}(F) = \mathbb{P} \left( \bigcap_{i = 1}^L F_{i} \right)  = \prod_{i=1}^L \mathbb{P}(F_{i} ~|~ F_{i-1}),
\label{eq:P_RE_decomp}
\end{equation}
with, by convention, $F_0 = \Theta$.  Computing $\Ptau$
according to~\eqref{eq:P_RE_decomp} requires an efficient and accurate method
for estimating the $L$ conditional probabilities.  We use
a modification of the Metropolis-Hastings algorithm to accomplish this~\cite{au2001estimation}. This modified Metropolis algorithm (MMA), which
belongs to the family of Markov Chain Monte Carlo (MCMC) methods, draws samples
from a conditional distribution and either accepts or rejects the samples based
on a chosen acceptance parameter. One uses samples that belong to $F_{i-1}$ as
seeds for estimating the conditional probability $\mathbb{P}(F_i | F_{i-1})$.
We refer the interested reader to \cite{beckzuev}, which provides a high level
discussion of MMA as well as other variants of SS; a more thorough analysis of
MMA and MCMC algorithms can be found in~\cite{papaioannou2015mcmc}.

Choosing a proper sequence of thresholds $\{\tau_i\}_{i = 1}^L$ is a major challenge of the SS method. Since one has little prior knowledge of the PDF of
$q(\boldsymbol\theta)$, it is often not feasible to prescribe the sequence of
thresholds a priori. Instead, one may require that $\mathbb{P}(F_{i}
~|~ F_{i-1}) = p_0$, $i = 1, \dots, L-1$, for some  chosen
quantile probability $p_0$~\cite{au2001estimation}.
We can then iteratively estimate the proper threshold at each ``level'' of the
algorithm; the SS estimator of \eqref{eq:P_RE_decomp}  takes the form 
\begin{equation} 
\Ptau \approx  \hat{P}_{\bar{\tau}}^{SS} = p_0^{L-1} ~\mathbb{P}(F_L ~|~ F_{L-1}),
\label{eq:Prob_RE_practical}
\end{equation}
where the final conditional probability $\mathbb{P}(F_L ~|~ F_{L-1})$ is
estimated via the MMA procedure mentioned earlier. 
Although $p_0 =
0.1$ is a standard choice in engineering applications, there has been significant work
done to determine optimal values for $p_0$; this, in general,  depends on the
 QoI under consideration. It has been shown that, for practical purposes, the optimal $p_0$ lies in
the interval $[0.1, 0.3]$ and that, within this interval, the efficiency of SS
is insensitive to the particular choice of $p_0$~\cite{zuev2012bayesian}. 
With the approach for computing the sequence of
thresholds in~\eqref{eq:P_RE_decomp}, each $\tau_i$ is a random variable, estimated via a finite number of conditional samples. Consequently the number of levels or iterations necessary to terminate SS is also random. For a sufficiently large number of samples,
the number of levels is given in~\cite{vsehic2020estimation} as 
\begin{equation}
\label{eq:SS_levelsneeded} L - 1 = \left\lfloor
\frac{\log{\Ptau}}{\log{p_0}} \right\rfloor.  
\end{equation}

\paragraph*{Implementation of subset simulation}

For completeness, we provide an algorithm outline for the SS method 
in Algorithm~\ref{alg:SS}.
We assume Gaussian inputs in the examples considered in this article; the SS method can 
however be applied to
non-Gaussian input distributions, see Appendix B
of~\cite{melchers2018structural} for details.
Additional information on the implementation of the SS algorithm, including the
MMA implementation, is for instance available in~\cite{au2001estimation,
vsehic2020estimation, ullmann2015multilevel}. As this MCMC
implementation reuses the input parameters from each previous level to estimate
the threshold for the next level, this method does not
require any burn-in samples to draw from the conditional distribution; it begins by sampling from the previous rare event domain.  On the theoretical side, the SS algorithm is asymptotically unbiased and $\hat{P}_{\bar{\tau}}^{SS}$ converges almost surely to the true rare event probability, $\Ptau$. For a detailed convergence analysis of SS and derivation of its statistical properties, see~\cite{au2001estimation}.

\renewcommand{\algorithmicrequire}{\textbf{Input:}}
\renewcommand{\algorithmicensure}{\textbf{Output:}}
\begin{algorithm}[H]
\small
\caption{Subset Simulation}
\label{alg:SS}
\begin{algorithmic}[1]
\REQUIRE Rare event threshold, $\bar{\tau}$, samples per level, $N_{SS}$, quantile probability, $p_0$, routine that evaluates QoI, $q(\boldsymbol{\theta})$
\ENSURE Estimate of rare event probability: $\hat{P}_{\bar{\tau}}^{SS}$
\STATE Draw $N_{SS}$ samples of $\boldsymbol{\theta}$ from appropriate distribution
\STATE Compute $N_{SS}$ samples of the QoI, compute $\tau_1$ as the $p_0$ quantile
\STATE Save the $\lfloor N_{SS} \cdot p_0 \rfloor$ inputs such that $q(\boldsymbol{\theta}) > \tau_1$ as seeds for the next level
\STATE $i \leftarrow 1$ \hfill\COMMENT{$i$ indicates the current level}
\WHILE{$\tau_i < \bar{\tau}$}
\STATE $i \leftarrow i + 1$
\STATE Sample $\boldsymbol{\theta}$ by creating $\lfloor N_{SS} \cdot p_0 \rfloor$ Markov Chains, each with length $\lfloor p_0^{-1} \rfloor$ \hfill\COMMENT{For details,~\cite{au2001estimation}}
\STATE Using MCMC seeds, evaluate the QoI and compute $\tau_i$ as the $p_0$ quantile
\STATE Save the $\lfloor N_{SS} \cdot p_0 \rfloor$ inputs such that $q(\boldsymbol{\theta}) > \tau_i$ as seeds for the next level
\ENDWHILE 
\STATE Using seeds from $F_{L-1}$, sample the QoI and estimate $\mathbb{P}(F_L ~|~ F_{L-1})$ using MC
\STATE Evaluate $\hat{P}_{\bar{\tau}}^{SS} = p_0^i ~\mathbb{P}(F_L ~|~ F_{L-1})$ \\
\end{algorithmic}
\end{algorithm}

\paragraph*{Computational cost}
We turn now to the computational cost of estimating $\Ptau$ using SS.
The computational cost is measured in terms of the number of function
evaluations required to run the algorithm. As the number of levels $L$ is random, so is the computational cost
associated with SS. For simplicity, we assume for our cost 
analysis that a sufficient number
of samples has been used so that $L$ does not vary. The total number of QoI
evaluations required by SS is $L \cdot N_{SS}$, where $N_{SS}$ is a
user-defined parameter that determines the number of samples per intermediate
level of the iteration. Say, for example, the true rare event probability is
$10^{-6}$ and we wish to estimate $\Ptau$ with a coefficient of
variation within $\delta = 0.1$. For standard MC sampling, we would need
$N \geq 1/(\delta^2 \cdot \Ptau) = 10^8$ samples of the QoI. Take
the SS method with a quantile probability of $p_0 = 0.1$. Then, according to
\eqref{eq:SS_levelsneeded}, we would have $L = 7$, corresponding to 7 levels of
conditional probabilities. The coefficient of variation for each of the
conditional probabilities is more difficult to quantify, however, as in the
case of the standard MC estimator, 
they are proportional to $1 / p_0$; see~\cite{au2001estimation}.
In this case, one would expect to see a significant reduction in the cost of estimating $\Ptau$ with SS.  

We lastly emphasize the power of SS for estimating rare event probabilities in
the context of QoIs with high-dimensional inputs. 
Not only does SS improve upon the slow
convergence rates of standard MC by a wide margin, it also inherits the property
of having a convergence rate independent of input dimension. 

\renewcommand{\SS}{\mathrm{SS}}

\section{Surrogates for GSA of rare event probabilities}\label{sec:SurrogatesforGSA}
We seek to apply  variance-based  GSA to  $\Ptau(\boldsymbol\xi)$, defined
in~\eqref{eq:qoixi}, with respect to components of $\boldsymbol\xi$.
To mitigate the computational expense of performing such analysis, we combine the SS
algorithm and surrogate models, in the form of 
polynomial chaos expansions (PCEs). We assume $\boldsymbol\xi$ to be an $M$-dimensional vector 
with independent entries. The procedure, which amounts to double-loop sampling, is  outlined below:  
	\begin{itemize}
\item Generate hyper-parameter samples 
$\{ \boldsymbol\xi^{(j)} \}_{j=1}^{\Ns}$
\item For each $j \in \{1, \ldots, \Ns\}$, estimate 
$\Ptau(\boldsymbol\xi^{(j)})$ using SS; denote 
these estimates by $\Ptau^{(j)} = \SS(\Ptau(\boldsymbol\xi^{(j)}))$
\item  Use the (noisy) 
function evaluations $\{\Ptau^{(j)}\}_{j=1}^{\Ns}$ 
to compute a surrogate model 
$\tilde{P}_{\bar\tau}(\boldsymbol\xi) \approx
\Ptau(\boldsymbol\xi)$
\item Compute the Sobol' indices of $\tilde{P}_{\bar\tau}(\boldsymbol\xi)$.
\end{itemize}
	
Instead of using SS for computing  $\Ptau(\boldsymbol\xi^{(j)})$, 
one may be tempted to apply a surrogate further ``upstream" by computing
 a surrogate model $\tilde{q}_{\boldsymbol\xi^{(j)}}(\boldsymbol\theta)$
for $q(\boldsymbol\theta)$ from samples 
$\{ q(\boldsymbol\theta^{(k)})
\}_{k = 1}^n$ drawn from law of 
$\boldsymbol\theta$ as determined by $\boldsymbol\xi^{(j)}$. 
This surrogate model of $q$ can then be used for fast approximation 
of the rare event probability $\Ptau(\xi^{(j)})$.
This procedure,
however, has two major pitfalls: (i)  an expensive surrogate modeling
procedure must be carried out for each $j \in \{1, \ldots, \Ns\}$
and, more importantly, (ii) surrogate models are 
typically poorly suited to the task of rare event estimation. Indeed, 
surrogates
typically fail to
capture the tail behavior of the distribution of the QoI $q$, making them
unsuitable for rare event simulations.
This shortcoming  is well-documented in the 
uncertainty quantification literature
~\cite{peherstorfer2017combining, li2010evaluation} although efforts are being made to tailor the
surrogate model construction process for the efficient estimation of
rare event probabilities~\cite{li2011efficient,li2010evaluation}.

\paragraph*{PCE surrogate for rare event probability}

Our approach  leverages the properties of PCE surrogates for fast estimation of
Sobol' indices~\cite{le2010spectral,crestaux2009polynomial}; it also  takes advantage of the regularity of the mapping $\boldsymbol\xi
\mapsto \Ptau(\boldsymbol\xi)$.  Specifically,  assuming the PDF of
$\boldsymbol\xi$ satisfies certain (mild) differentiability and integrability
conditions, one can show that $\Ptau(\boldsymbol\xi)$ is a
differentiable function of $\boldsymbol\xi$; see~\cite[Proposition
3.5]{AsmussenGlynn07}.

The PCE of $\Ptau(\boldsymbol\xi)$ is defined as
\begin{equation}
\label{eq:PCE_def} \tilde{P}_{\bar\tau}(\boldsymbol{\xi}) = \sum_{k = 0}^{\Npc} \beta_k
\Psi_k(\boldsymbol{\xi}), 
\end{equation} 
where $\Psi_0, \cdots, \Psi_{\Npc}$ belong to a
family of orthogonal polynomials and $\beta_0, \cdots, \beta_{\Npc}$ are the 
(scalar) PCE coefficients. The specific family of polynomials is chosen
to guarantee orthogonality with respect to the PDF of $\boldsymbol{\xi}$; 
see, e.g.,~\cite{le2010spectral}. We use a total order truncation 
scheme for the PCE: the multivariate polynomial basis contains all possible
polynomial basis elements up to a total polynomial order $r$. In this 
case, $\Npc$ in~\eqref{eq:PCE_def} satisfies
\[ 
\Npc + 1 = \frac{(M+r)!}{M! r!}. 
\] 
The coefficients $\beta_0, \dots,
\beta_{\Npc}$ can be computed in a number of ways, including 
non-intrusive spectral projection or 
regression~\cite{le2010spectral,BlatmanSudret10,crestaux2009polynomial}.
A regression based approach is preferred here because 
the evaluations of $\Ptau$ are noisy due to sampling errors incurred in the SS 
procedure. We estimate the vector $\boldsymbol{\beta} = 
[\begin{matrix} \beta_0, \, \beta_1, \, \cdots, \beta_{\Npc}\end{matrix}]$
from function evaluations
$\Ptau^{(j)} = \SS(\Ptau(\boldsymbol\xi^{(j)}))$, $j = 1, \ldots, 
\Ns$,
by solving the penalized least squares problem
\begin{equation}
\label{eq:PCE_sparse} 
\min_{\boldsymbol{\beta}}
\sum_{j=1}^{\Ns} 
\big[
\Ptau^{(j)} - \sum_{k = 0}^{\Npc} \beta_k \Psi_k(\boldsymbol{\xi}^{(j)})
\big]^2 \quad \mathrm{s.t.} \quad ||\boldsymbol{\beta} ||_1 \leq \lambda.
\end{equation}
In (\ref{eq:PCE_sparse}), the penalty parameter  $\lambda$
acts as a sparsity control on the recovered PCE coefficients. We generate  the realizations $\{ 
\boldsymbol{\xi}^{(j)}\}_{j=1}^{\Ns}$ of
the hyper-parameter vector through  Latin
hypercube sampling;
for further details on the implementation of sparse
regression for PCE, see~\cite{fajraoui2017sequential, hampton2016compressive}. The numerical results in Section~\ref{sec:NumericalResults} are obtained using the SPGL1 solver~\cite{spgl1site}.

\paragraph*{GSA of $\Ptau$ using the PCE surrogate}
As is well-known, the Sobol' indices of a PCE surrogate can be computed analytically.
For example, the first order Sobol' indices, $S_i(\Ptau)$, $i = 1, \ldots, M$, 
of $\Ptau$ can be approximated as follows: 
\begin{equation}\label{eq:Sobol_PCE}
S_i(\Ptau) \approx S_i(\tilde{P}_{\bar\tau}) = \frac{\sum_{k \in
K_i} \beta_k^2~\mathbb{E}[\Psi_k^2] }{\sum_{k = 1}^{\Npc} \beta_k^2
~\mathbb{E}[\Psi_k^2] }, 
\end{equation}
where $K_i$ denotes the set of all PCE terms that depend only on $\xi_i$. Sobol'
indices for arbitrary subsets of variables, as well as total indices, can be
computed in an analogous manner~\cite{le2010spectral,Alexanderian13}.  
In practice, PCE surrogates with modest accuracy are often sufficient
to obtain reliable estimates of Sobol' indices. 

While the above approach for GSA of $\Ptau$ does require repeated
simulations of the QoI $q$ during the calls to the SS algorithm, it still provides
orders of magnitude speedup over the standard ``pick and freeze'' MC methods, also known as Saltelli sampling,
for computing the Sobol' indices of $\Ptau$~\cite{saltelli2010}.  Indeed,  a fixed sample  $\{
\boldsymbol{\xi}^{(j)}\}_{j=1}^{\Ns}$ with modest $\Ns$ is sufficient to compute the
PCE surrogate from which the Sobol' indices can be computed at a negligible
computational cost. Moreover, the sparse regression approach for estimating PCE
coefficients is forgiving of noisy function evaluations. Therefore, large sample sizes are not needed in the calls to the SS algorithm. 
We demonstrate the merits of the proposed approach in our computational 
results presented in Section~\ref{sec:NumericalResults}. 





\section{Numerical results}\label{sec:NumericalResults}
We summarize, in Section~\ref{sec:ResultsAnalytic}, the computational results
for the motivating example from Section~\ref{sec:motivatingexample}; a more
challenging model problem involving flow through porous media is considered in
Section~\ref{sec:SubsurfaceFlow}.

\subsection{Results for the analytic test problem}\label{sec:ResultsAnalytic}

\begin{figure}[!htb]
\centering
\includegraphics[scale=.5]{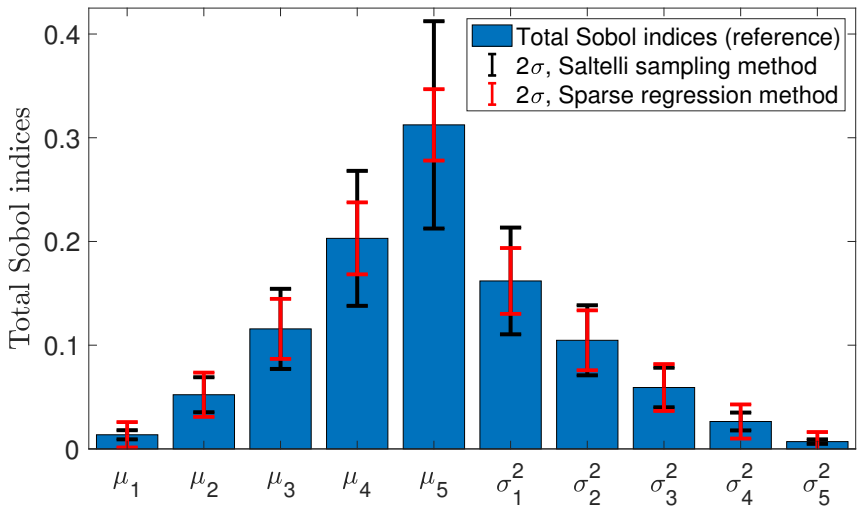}
\caption{Total Sobol' indices of $\Ptau$, with $\bar\tau=3$,
from (\ref{eq:exactp}); the error bars illustrate the variability of the two sampling methods (Saltelli sampling and sparse regression PCE) around the reference values.}
\label{fig:Sobol_errorbar}
\end{figure}

We consider the example from Section~\ref{sec:ResultsAnalytic} and 
study $\Ptau$ with $\bar\tau = 3$.
To establish a baseline for the values of the Sobol' indices of
$\Ptau(\boldsymbol\xi)$, we compute the total order Sobol'
indices directly from~\eqref{eq:exactp} using Saltelli sampling. The reference
Sobol' indices are computed with $10^6$ samples for each of the conditional terms;
convergence was numerically verified. We plot the reference total indices in Figure~\ref{fig:Sobol_errorbar} for comparison. We now compare the reference indices with those obtained through the PCE surrogate when $\Ptau(\boldsymbol\xi)$ is computed analytically using Equation \eqref{eq:exactp}. 
We allocate $10^3$ samples of $\Ptau(\boldsymbol\xi)$ each for the Saltelli sampling method and sparse regression PCE method. The Saltelli method requires $N(d+1)$ samples~\cite{saltelli2010} and so we divide the budget of $10^3$ samples equally among each conditional term. Each PCE coefficient can be estimated using the full set of $10^3$ samples. For a fair comparison, we use Latin hypercube sampling for both the PCE and Saltelli method. We also use a total PCE order of 3 and the penalty parameter $\lambda = 5\times 10^{-2}$.  Given that the set of total indices is computed, in each method, using a finite
number of samples, each index is a random variable with an associated
distribution. We compare the standard deviation of each total index for the two
GSA methods. In each case, we compute $10^3$ realizations of the full set of
total indices and compare their respective standard deviations in
Figure~\ref{fig:Sobol_errorbar}. 

Figure \ref{fig:Sobol_errorbar} illustrates the higher accuracy, or lower
variance, of PCE with sparse regression over Saltelli sampling: the standard deviation of the largest Sobol'
index is roughly 3 times smaller with sparse regression than it is with Saltelli sampling. 
 This gap in accuracy
appears to diminish for smaller indices, although the methods do not show comparable accuracy until
the indices are below 0.1. As
$\Ptau$ can be expressed analytically, there may be  additional benefits of the
sparse regression method to be seen when one considers performing GSA on a rare
event probability with noise due to sampling. We note that the total order
of the PCE basis and the penalty parameter $\lambda$, which are user-defined
parameters, can be changed without the need for additional runs of SS. These
parameters can be cross validated in a post-processing step after the rare event
simulation step, providing flexibility in this approach without adding any
significant computational burden.

\begin{figure}[ht!]
\centering
\includegraphics[width=.55\textwidth]{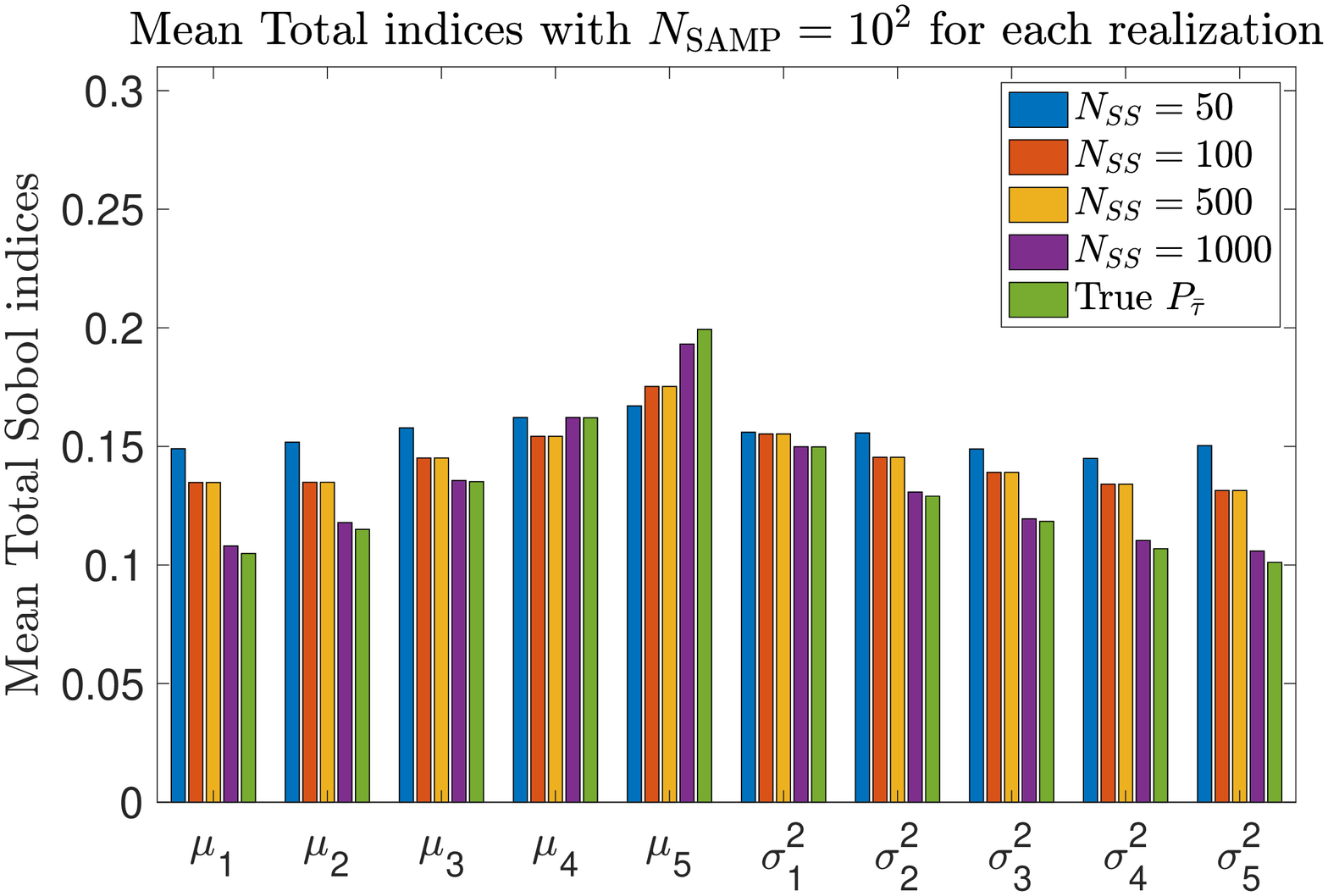} \\ 
\includegraphics[width=.55\textwidth]{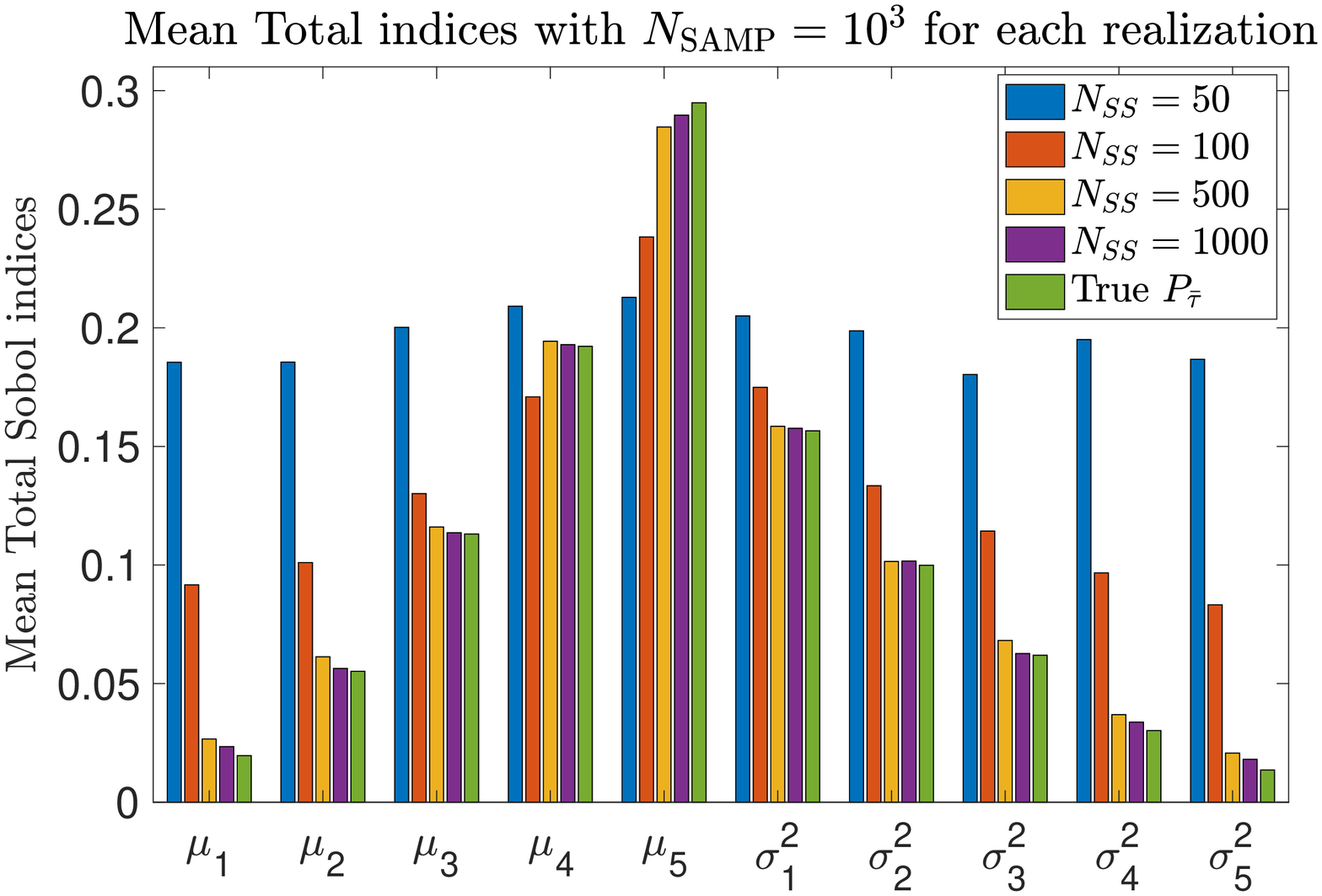} \\ 
\includegraphics[width=.55\textwidth]{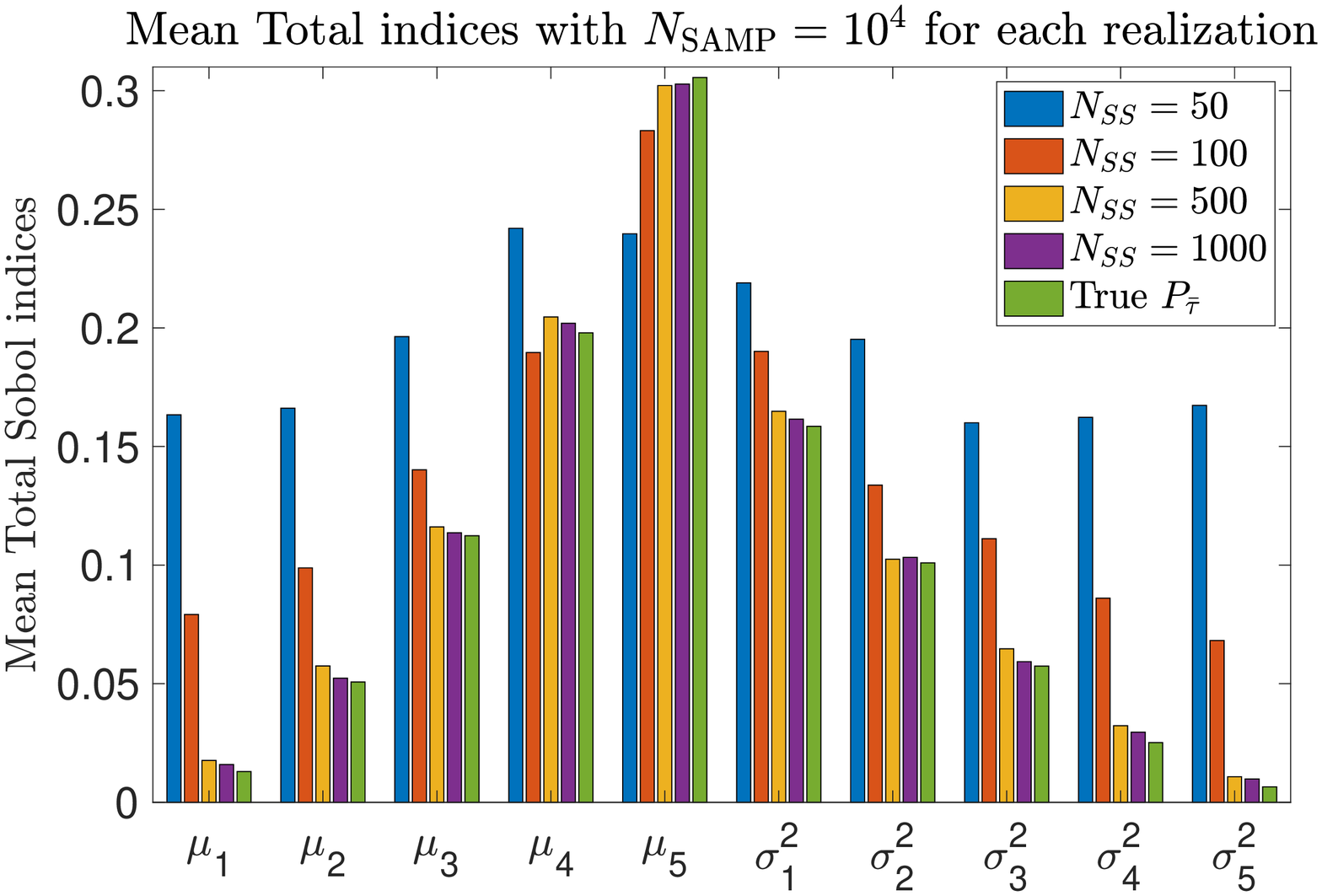}  
\caption{Mean Total Sobol' indices, varying the computational cost of SS and the PCE construction. Each plot varies $N_{\mathrm{SAMP}}$ and each colored bar varies $N_{SS}$, with the final bar of each index corresponding to the analytic $P_{\bar\tau}$.}
\label{fig:TSIs_compareNSS}
\end{figure}

When combining PCE-based GSA with SS for estimating $\Ptau(\boldsymbol\xi)$, there is a tradeoff 
between the inner loop cost of estimating $P_{\bar\tau}$ via SS and the outer
loop of aggregating $P_{\bar\tau}$ samples to build the PCE. In Figure
\ref{fig:TSIs_compareNSS}, we separately vary $N_{\mathrm{SS}}$ and
$N_{\mathrm{SAMP}}$ and examine the resulting distribution of the total Sobol'
indices, computed via sparse regression PCE. For a fixed $N_{\mathrm{SAMP}}$,
we compute multiple realizations of the total indices for several values of 
$N_{\mathrm{SS}}$.
Figure \ref{fig:TSIs_compareNSS}~(top) displays the expected value of the
total indices for $N_{\mathrm{SAMP}} = 100$. Regardless of how
accurately we estimate $\Ptau$, the indices do not approach their true values because 
the PCE is built using an inadequate number of samples, resulting
in a poor surrogate. By contrast, 
Figure~\ref{fig:TSIs_compareNSS}~(middle) shows that for $N_{\mathrm{SAMP}} = 10^3$,
we only need a modest $N_{\mathrm{SS}}$ to approximate the Sobol' indices.
Indeed, for $N_{\mathrm{SS}} = 500$, we are able to resolve the total indices
very well. We also examine the case of $N_{\mathrm{SAMP}} = 10^4$ in
Figure~\ref{fig:TSIs_compareNSS}~(bottom). Again, we are  able to resolve
the total indices well using only $N_{\mathrm{SS}} = 500$ and are able to
achieve the correct ordering for as little as $N_{\mathrm{SS}} = 100$.

These results indicate that 
(i)  a modest number of samples allocated to SS is enough 
to get a rough estimate of $\Ptau$ and (ii)  a moderate number of realizations
of $\Ptau(\boldsymbol{\xi})$ is then sufficient for accurate GSA. 
In other
words, given rather poor estimations of $\Ptau$, we are still able to extract
accurate GSA results, due to the fact that the sparse regression technique is robust to noisy QoI evaluations.

\subsection{Subsurface flow application}\label{sec:SubsurfaceFlow}
We consider the equations for single-phase, steady state flow in a square domain $\mathcal{D} = [0,1]^2$:
\begin{equation}\label{equ:Poisson}
\begin{aligned}
-\nabla \cdotp \left( \frac{\kappa}{\mu} \nabla p \right) &= 0 \quad \text{in} ~\mathcal{D}, \\
p & = 1 \quad \text{on} ~\Gamma_{1}, \\
p & = 0 \quad \text{on} ~\Gamma_{2}, \\
\nabla p \cdotp n &= 0 \quad \text{on} ~\Gamma_{3}, 
\end{aligned}
\end{equation}
where $\kappa$ is the permeability, $\mu$ is the viscosity, and $p$ is the
pressure. The boundaries $\Gamma_1, \Gamma_2$, and $\Gamma_3$ indicate the left
boundary, the right boundary, and the top/bottom boundaries, respectively. The
Darcy velocity is defined as $\boldsymbol{v} = -\frac{\kappa}{\mu} \nabla p$. 
In the present study, we let
$\mu = 1$. The source of uncertainty in this problem is
in the permeability field, which we model as a random field.  We consider the
flow of particles through the medium and focus on  determining the probability
that said particles do not reach the outflow boundary in a given amount of
time. 
This problem has been used previously as a
test problem for rare event estimation in~\cite{ullmann2015multilevel} as it
pertains to the long-term reliability of nuclear waste repositories.  Our goal
is to perform GSA with respect to the hyper-parameters that define the
distribution law of the permeability field.

\paragraph*{The statistical model for the permeability field}
Following standard practice~\cite{ullmann2015multilevel, cleaves2019derivative}, we model the permeability field as a log-Gaussian random field: 
\begin{equation} \label{eq:logkappa}
\log{\kappa}(x, \omega) = a(x, \omega) = \bar{a}(x) + \sigma_a z(x, \omega),
\end{equation}
where $x \in \mathcal{D}$ and $\omega$ belongs to 
sample space that carries the random process. Here, $\bar{a}$ is the 
mean of the random field, $\sigma_a$ is a scalar which controls the 
pointwise variance of the field, and $z$ is a centered (zero-mean) 
random process. 
We let the covariance function of $z$ be given by 
\begin{equation}
c_z(x, y) = \exp{\left(-\frac{|x_1 - y_1|}{\ell_x} - \frac{|x_2 - y_2|}{\ell_y}\right)}, \quad x, y \in \mathcal{D},
\label{eq:KLE_covariance}
\end{equation}
where $\ell_x$ and $\ell_y$ denote the correlation lengths in horizontal and
vertical directions.
The random field is represented via a 
truncated Karhunan-L\`{o}eve expansion (KLE):
\begin{equation}
a(x, \omega) \approx \bar{a}(x) + \sum_{k = 1}^{\NKL} 
\sqrt{\lambda_k} ~\theta_k(\omega) ~e_k(x).
\label{eq:KLE}
\end{equation}
In this representation, 
$\theta_1, \ldots, \theta_{\NKL}$ are independent standard normal
random variables and 
$(\lambda_i, e_k)$, $k = 1, \ldots, \NKL$, are the leading 
eigenpairs of the covariance operator of the process.
Our  setup for the uncertain log-permeability field follows 
the one
in~\cite{cleaves2019derivative}: we use permeability data
from the Society for Petroleum Engineers~\cite{spe_data} to define $\bar{a}$.
Once we truncate the KLE, the random vector
$\boldsymbol\theta = [\begin{matrix} 
\theta_1 & \theta_2 & \cdots & \theta_{\NKL}\end{matrix}]^\top$ fully
describes
the uncertainty in the log-permeability field. 

To ensure that the KLE accurately models the variability of the
infinite-dimensional field, we examine the eigenvalue decay of the covariance
operator with the goal to truncate the KLE so
that at least 90\% of the average variance of the field is maintained. 
For $\ell_x = \ell_y = 0.4$, which are the smallest correlation lengths considered in 
the present study, we require at least $\NKL = 126$. The number of retained KL
modes then determines the dimensionality of the rare event estimation
problem, and is henceforth fixed at $126$. The dimension independent properties of SS is advantageous in this
regime.

\begin{figure}[ht!]
\centering
\includegraphics[width=.35\textwidth]{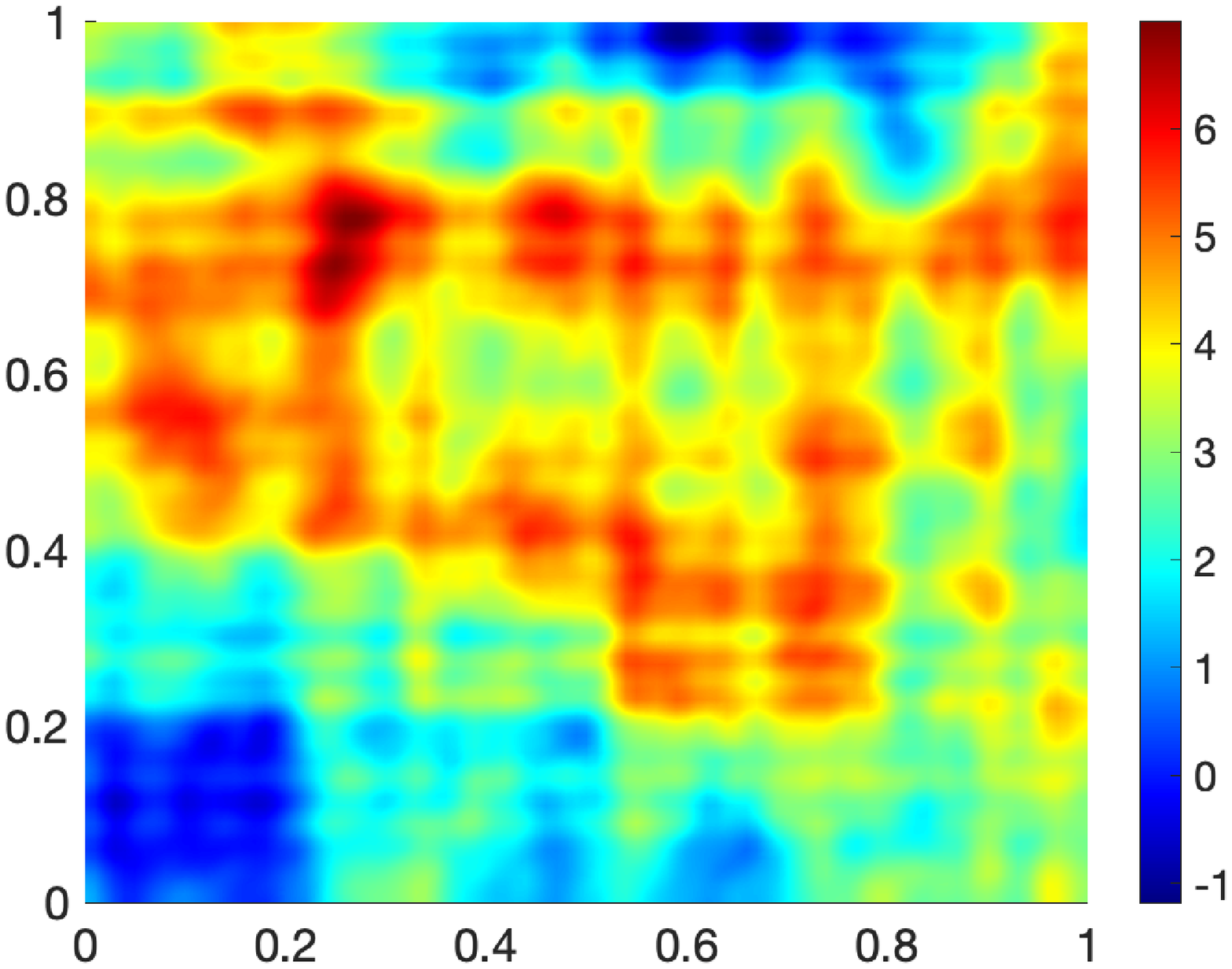}
\includegraphics[width=.35\textwidth]{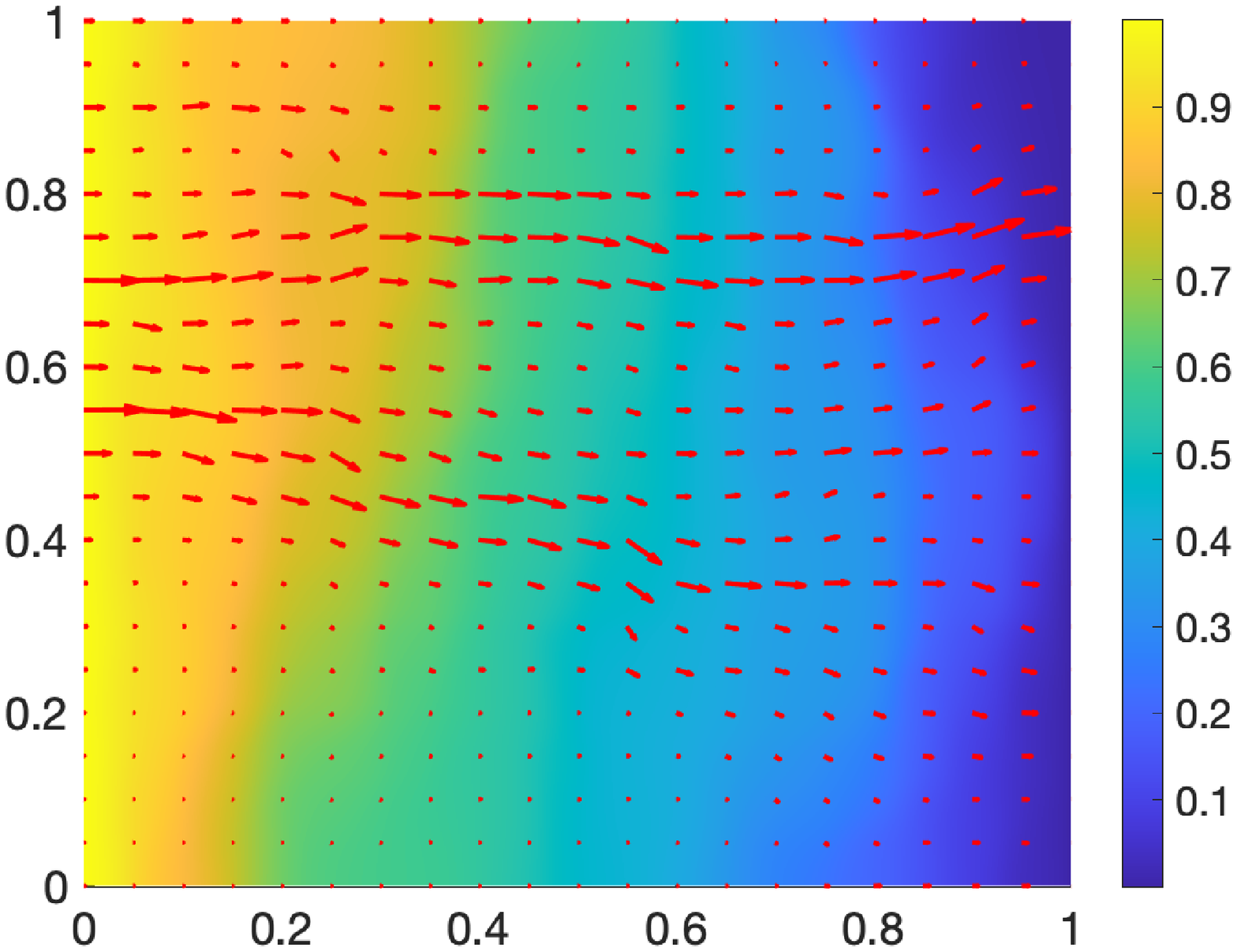}
\includegraphics[width=.35\textwidth]{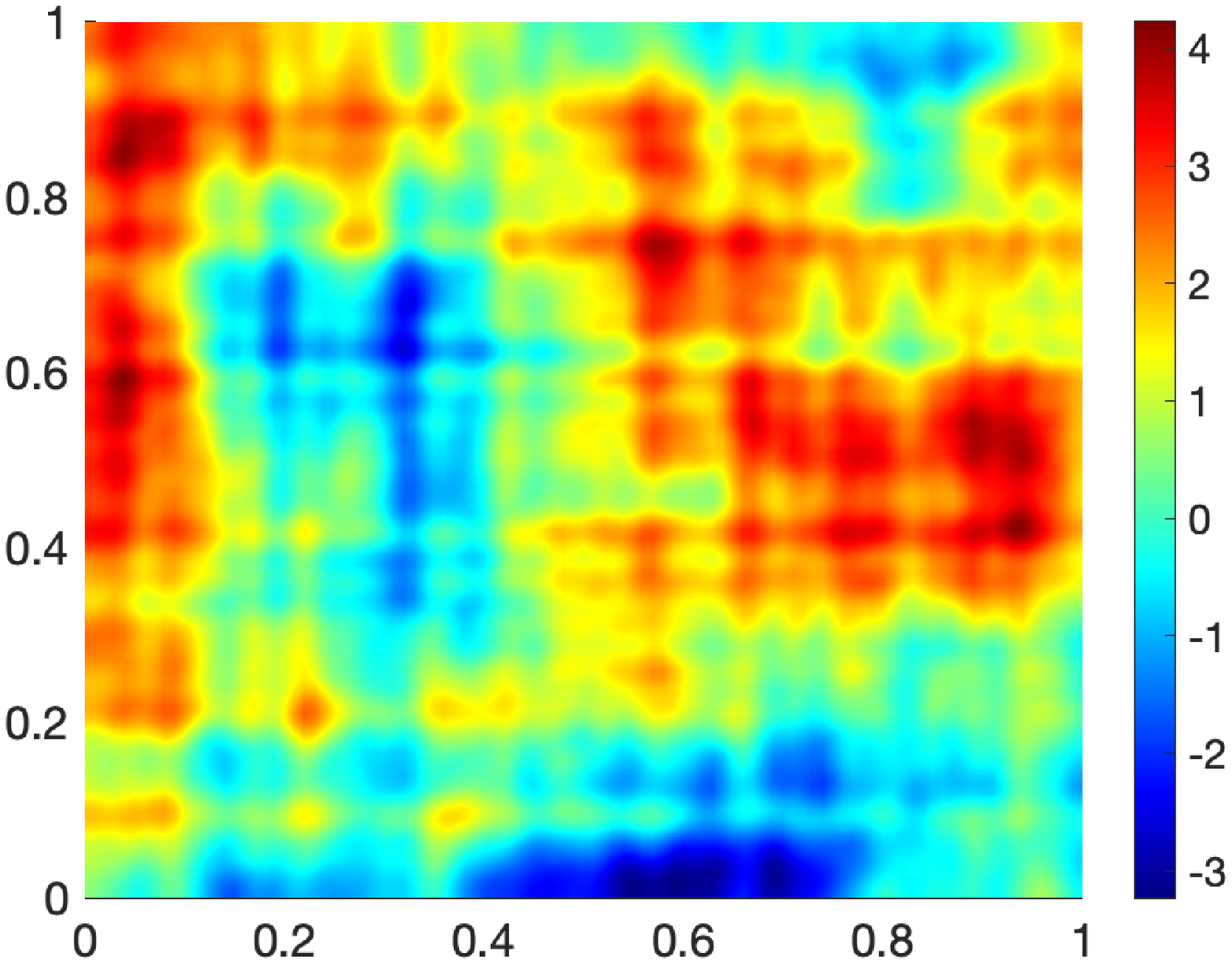}
\includegraphics[width=.35\textwidth]{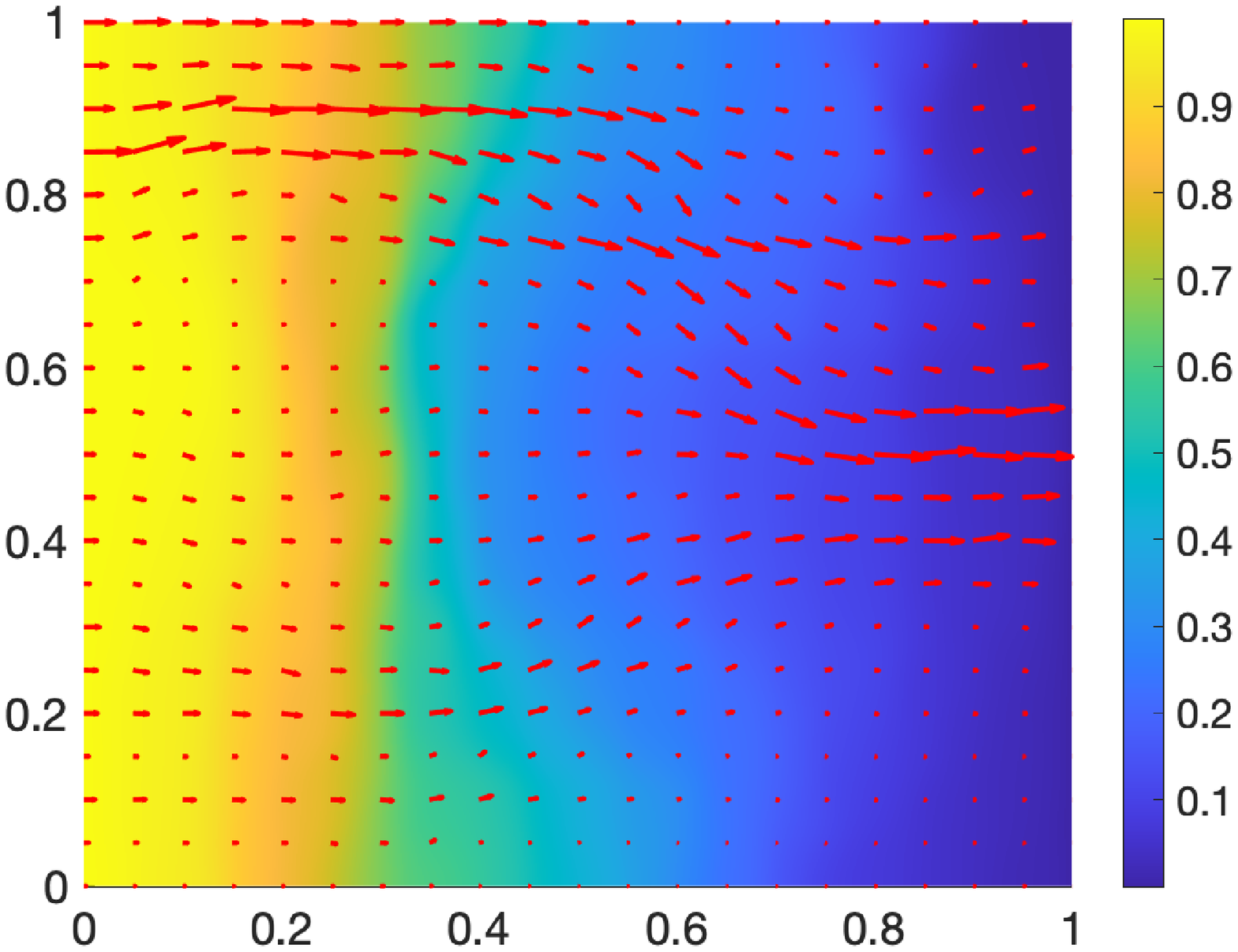}
\caption{Left: plots showing two realizations of the log permeability field. Right: the corresponding pressure solution and arrows indicating the resulting Darcy velocity field.}
\label{fig:KLE_Darcyvel}
\end{figure}

For illustration, we plot two realizations of the random field, 
with the corresponding pressure and velocity fields obtained
by solving the governing PDE~\eqref{equ:Poisson}, in Figure~\ref{fig:KLE_Darcyvel}.
In our computations, we solve the PDE 
using piecewise linear finite elements in
\textsc{Matlab}'s finite element toolbox with 50 mesh points in each direction.

\paragraph*{The QoI and rare events under study}
 The position $\boldsymbol{x}$
of a particle moving with the flow through the medium is determined by 
the following ODE
\begin{equation}
\begin{aligned}
\frac{d\boldsymbol{x}}{dt} &= \boldsymbol{v}, \\
\quad \boldsymbol{x}(0) &= \boldsymbol{x}_0,
\end{aligned}
\label{eq:ODE_darcyflow}
\end{equation}
where $\boldsymbol{v}$ is the Darcy velocity.
We consider a single particle with initial position 
at 
$\boldsymbol{x}_0 = 
\begin{bmatrix} 0 \\ 0.5\end{bmatrix}$. 
The solution $\boldsymbol{x}$ of~\eqref{eq:ODE_darcyflow} 
depends not only on time but also on $\boldsymbol\theta$ due to dependence of 
$\boldsymbol{v}$ on $\boldsymbol\theta$, i.e., 
$\boldsymbol{x} = \boldsymbol{x}(t, \boldsymbol\theta)$.
We take  the QoI $q$ as the hitting time, i.e.,  the time it takes a particle to travel through the medium from left to right
\[
q(\boldsymbol\theta) = \{ t : x_1(t, \boldsymbol\theta) = 1 \}.
\]
We aim to determine the rare event probability $P_{\bar{\tau}} = 
\mathbb{P}(q > \bar{\tau})$. The parameters $\ell_x$, $\ell_y$,  and $\sigma_a$ parametrize the uncertainty
in the permeability field; we consider them as hyper-parameters and set
$\boldsymbol{\xi} = [\begin{matrix} \ell_x & \ell_y &
\sigma_a\end{matrix}]^\top$.  We set  the nominal values of the
hyper-parameters $\boldsymbol{\xi}_{nom} = [\begin{matrix} 0.4 & 0.4 &
0.8\end{matrix}]^\top$. 
We simulate realizations of the permeability field at these nominal hyper-parameters and plot the distribution of $q$. Each of these realizations requires one PDE solve and one ODE solve.
\begin{figure}[ht!]
\centering
\includegraphics[width=.49\textwidth]{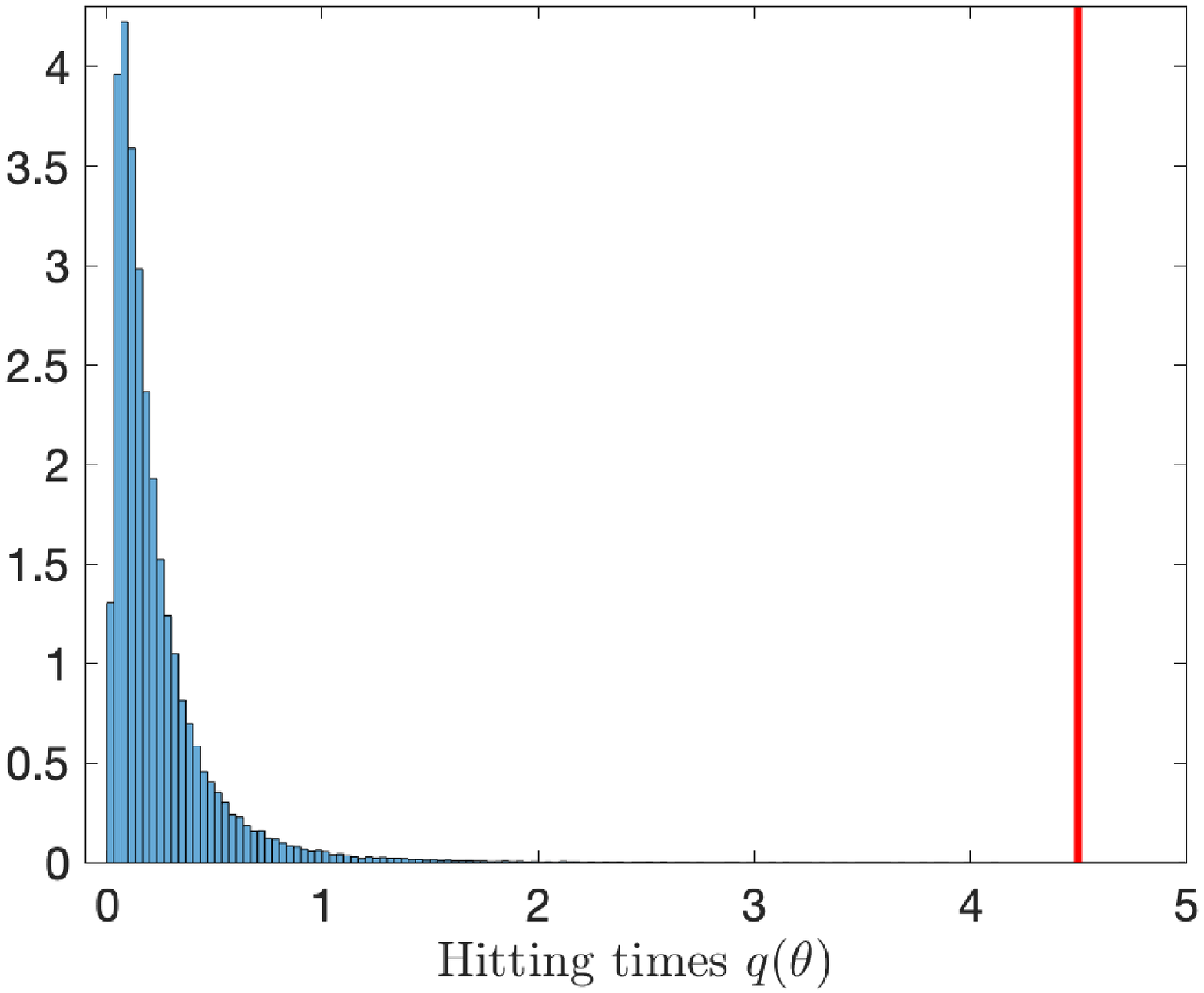}
\includegraphics[width=.49\textwidth]{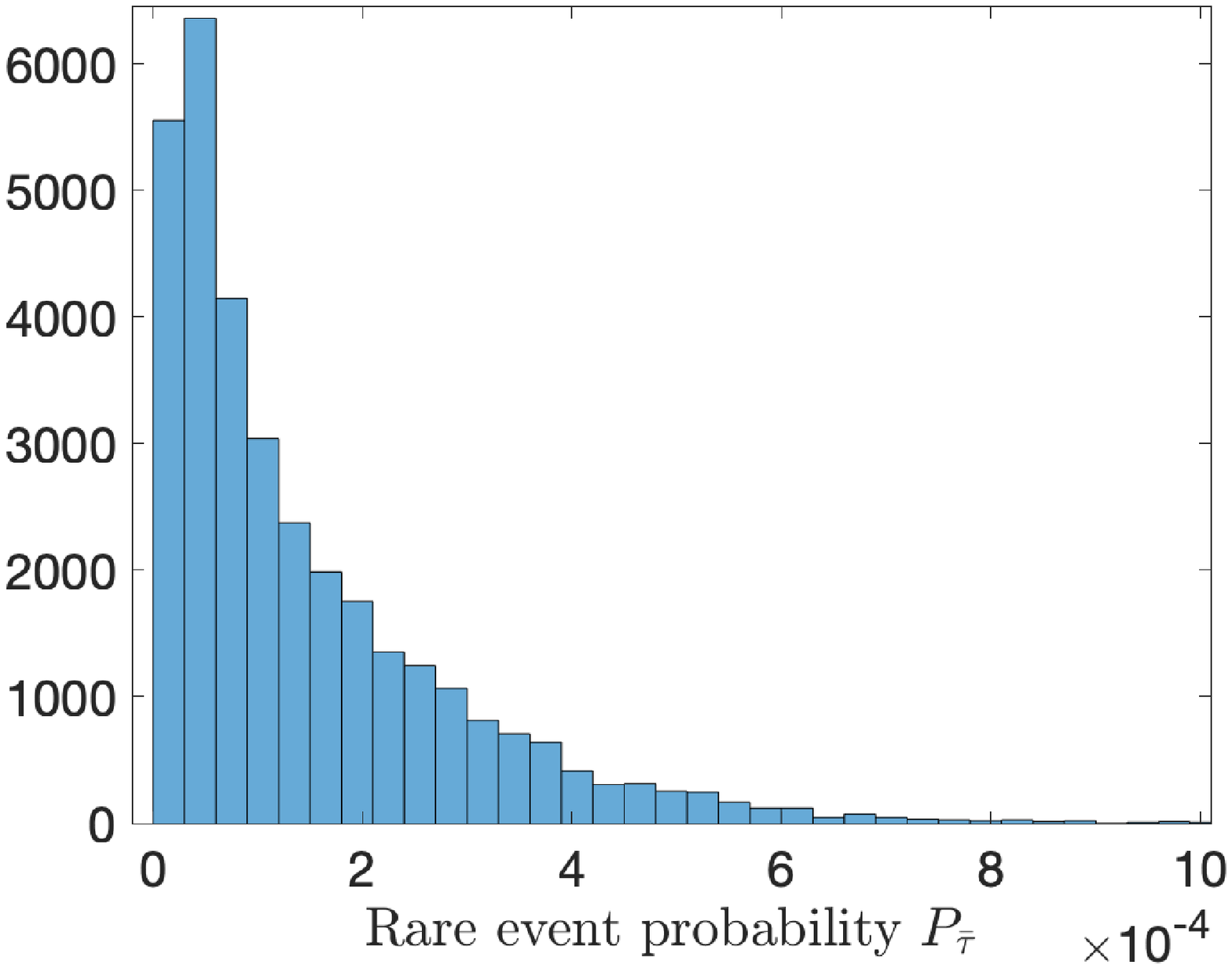}
\caption{Left: Histogram of $q$ for nominal hyper-parameters. Vertical line indicates rare event threshold of $\bar{\tau} = 4.5$. Right: histogram of the rare event probability, estimated via SS with uniformly distributed hyper-parameters.}
\label{fig:histogram_tauh}
\end{figure}
As illustrated in Figure~\ref{fig:histogram_tauh}, the distribution for $q$ corresponds to a heavy-tailed distribution. We select as the threshold $\bar{\tau} = 4.5$ and consider quantifying the sensitivity of $P_{\bar{\tau}}(\boldsymbol{\xi})$ with respect to the hyper-parameters defining the KLE. 

\paragraph*{Rare event probabilities and GSA}
In our first set of experiments, we use SS with $N_{SS} = 10^3$ samples per
intermediate level; each of these samples corresponds to one solution of the
full subsurface flow problem, including a PDE and ODE solve.  For each
evaluation of SS, approximately 5 intermediate levels are used, resulting in
approximately $5 \times 10^3$ function evaluations. Our hyper-parameters are
drawn from a uniform distributed centered at $\boldsymbol{\xi}_{nom}$ with a
spread of plus or minus 10\% of $\boldsymbol{\xi}_{nom}$. We use these SS
estimations of $P_{\bar{\tau}}(\boldsymbol{\xi})$ in order to build the
corresponding PCE surrogate, where the polynomial basis is truncated at a total
polynomial order of 5. Note the decision of where to truncate the PCE basis
does not need to be made prior to estimating the set of
$P_{\bar{\tau}}(\boldsymbol{\xi})$ samples.

The samples for the hyper-parameters are drawn using a Latin hypercube sampling
scheme.  We use $10^3$ estimations of $P_{\bar{\tau}}(\boldsymbol{\xi})$ to
construct the PCE surrogate. Again, we use sparse regression to recover the
PCE coefficients, while promoting sparsity in the set of PCE coefficients,
and so mitigating the effects of noise induced by SS. 
\begin{figure}[ht!]
\centering
\includegraphics[width=.49\textwidth]{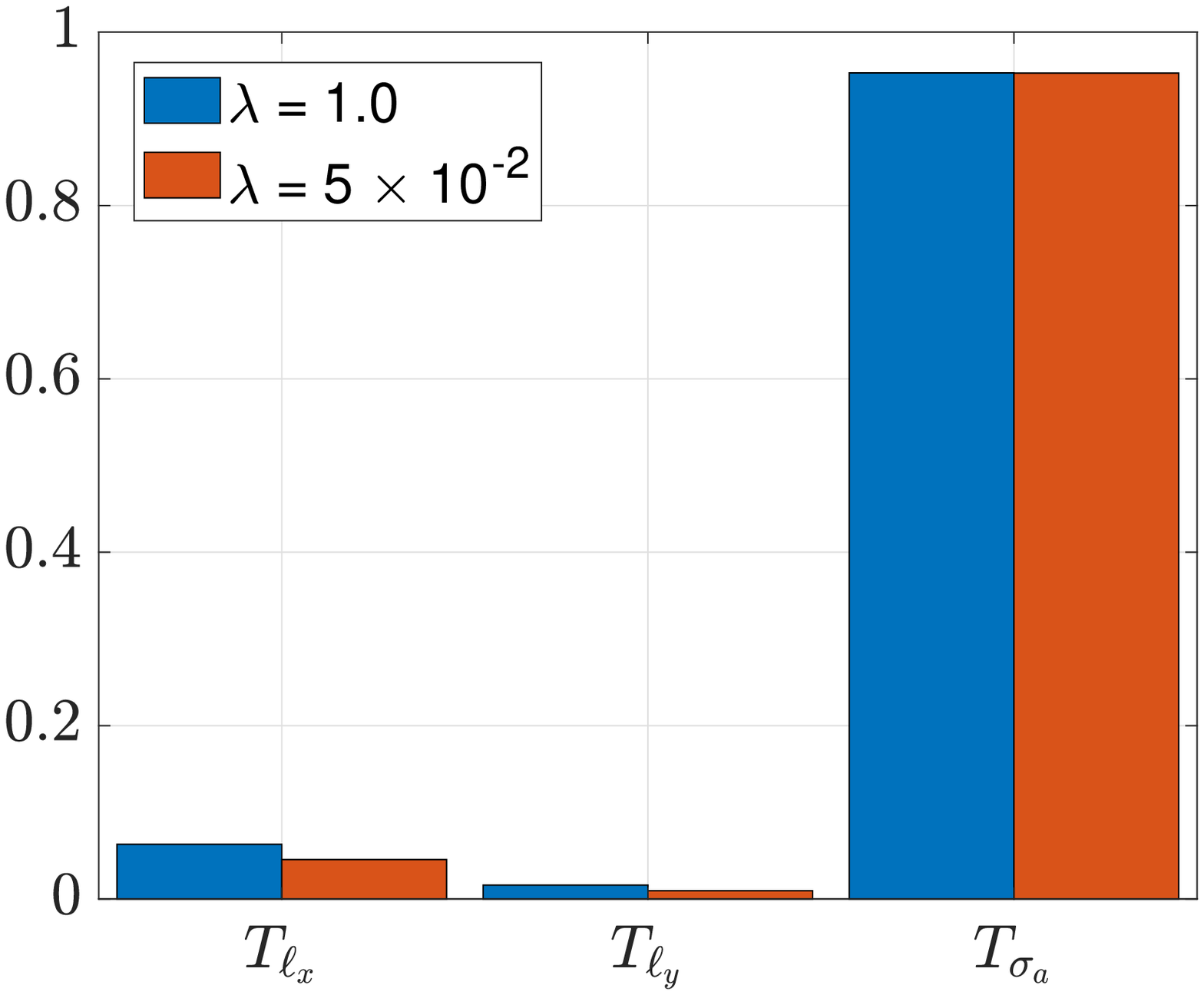}
\includegraphics[width=.49\textwidth]{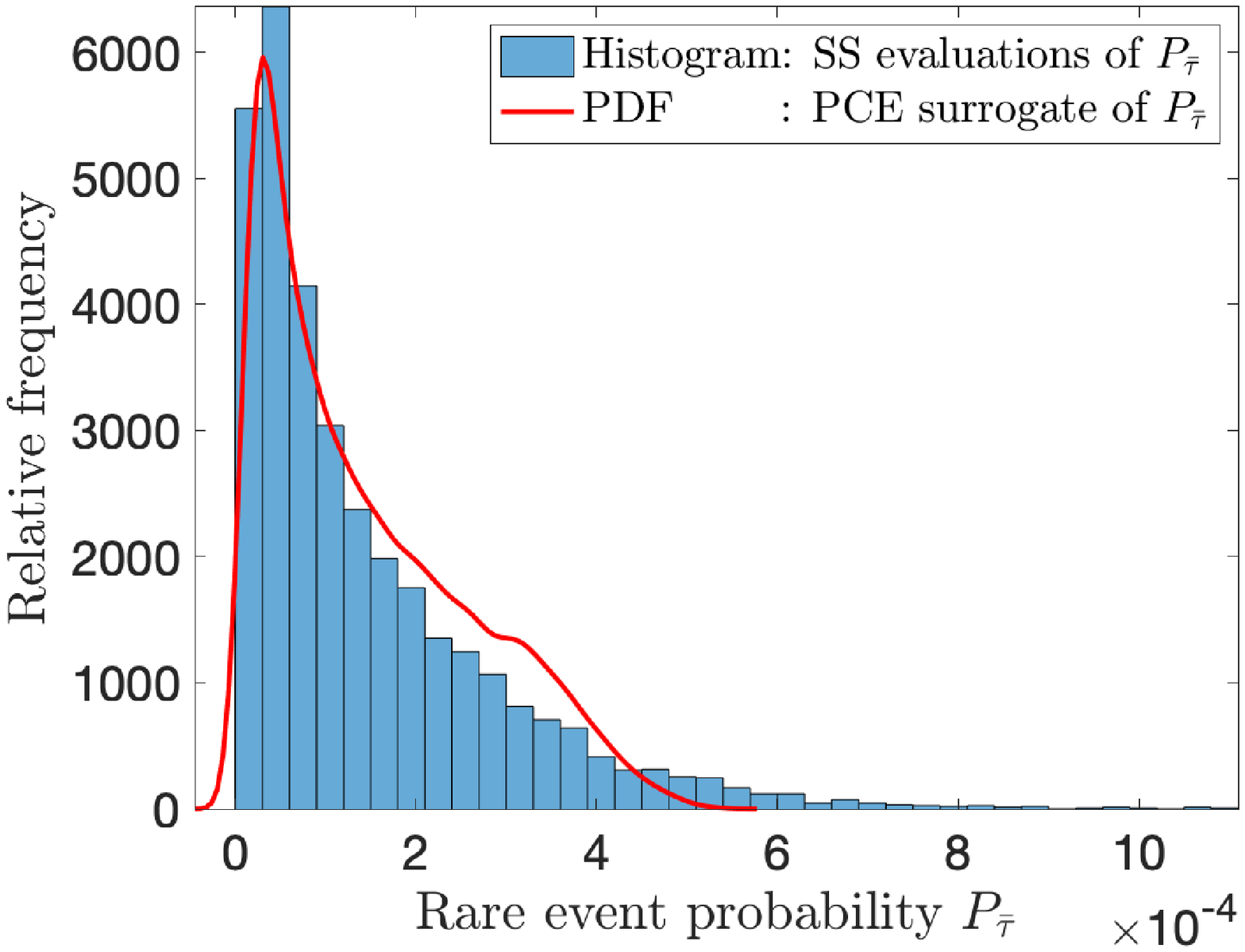}
\caption{
Left: total Sobol' indices for $P_{\bar{\tau}}(\boldsymbol{\xi})$ computed from recovered PCE coefficients; results are reported with regularization constant 
$\lambda = 1$ and $\lambda = 5\times10^{-2}$. Right: PDF of PCE surrogate compared with $\Ptau$ evaluation histogram. Used $\Ns = 10^4$ for better resolution of distributions.}
\label{fig:Sobol_Subsurf}
\end{figure}
In Figure \ref{fig:Sobol_Subsurf}, we use two different values of $\lambda$
when promoting sparsity in order to illustrate the effect of $\lambda$ on the results. 
Note that when $\lambda$ is made smaller, the PCE
coefficients decrease in magnitude, promoting a sparser PCE
spectrum. In both cases, the ordering of the total Sobol' indices remains
consistent, and thus, conclusions with respect to parameter sensitivity are
unaffected. For this experiment, we therefore conclude that choosing $\lambda$ by trial and error is sufficient. Should one encounter a scenario where the GSA results are more 
sensitive to $\lambda$, more systematic approaches are possible~\cite{BlatmanSudret10, hampton2016compressive}. 
\begin{figure}[ht!]
\centering
\includegraphics[width=.49\textwidth]{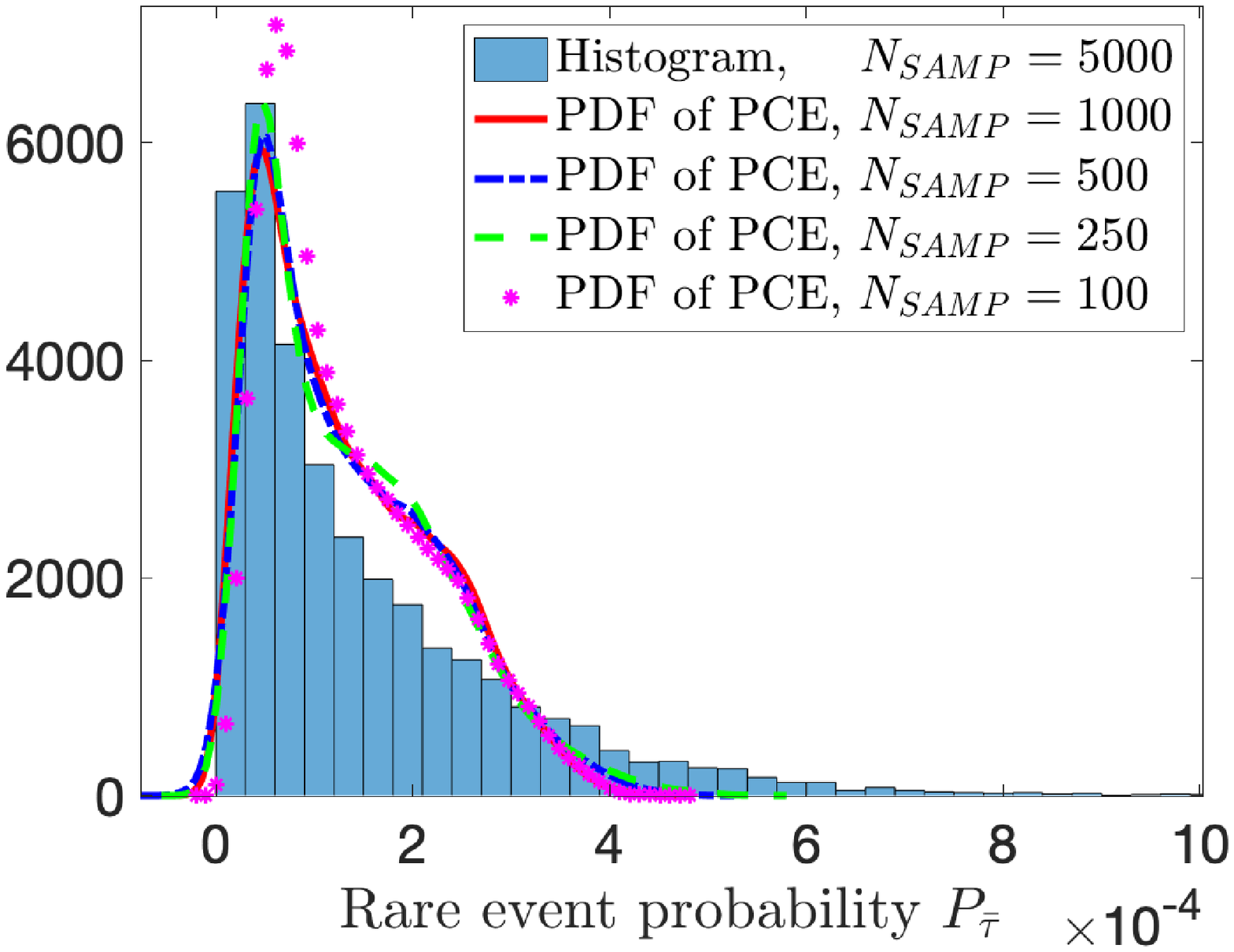}
\includegraphics[width=.49\textwidth]{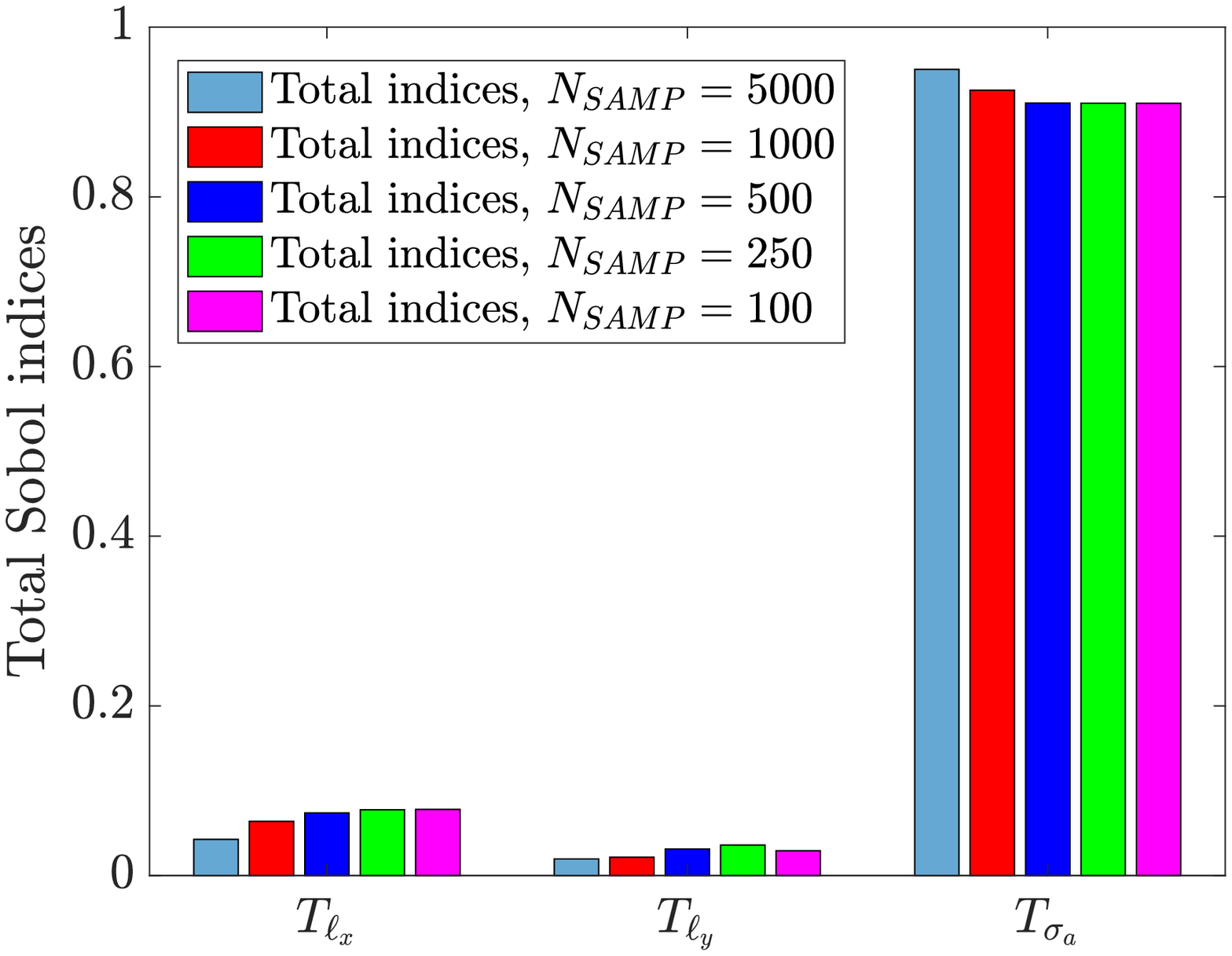}
\caption{Distributions of $\Ptau$ for $\NSS = 500$ and varying values of $\Ns$. For each $\Ns$, we build the PCE surrogate and approximate its PDF with $10^5$ samples. The set of $\Ptau$ samples used for differing $\Ns$ is nested within sets of larger samples. Corresponding total indices are included, computed directly from the PCE surrogates.}
\label{fig:PDF_compare_darcyflow}
\end{figure}

We lastly return to the key point made in Section \ref{sec:ResultsAnalytic},
that the proposed method is capable of producing reliable GSA results, while
using a modest number of inner and outer loop samples ($\NSS$ and $\Ns$,
respectively). 
In Figure~\ref{fig:PDF_compare_darcyflow}, we report results corresponding to $\NSS =
500$. In the left panel of the Figure we study the effect of $\Ns$ on the
PDF of the PCE
surrogate. In the right panel, we plot the Sobol' indices
corresponding to each of the computed surrogates. The results in 
Figure~\ref{fig:PDF_compare_darcyflow} should also be compared with 
those in Figure~\ref{fig:Sobol_Subsurf}, where larger values of $\NSS$ and
$\Ns$ were used.
This experiment indicates that 
$\Ptau$ and the Sobol' indices themselves can be well-approximated with a
modest number of samples in both the inner and outer loops. In this case, using
both $\NSS$ and $\Ns$ on the order of $10^2$ is sufficient for obtaining
accurate GSA results. The combined cost of this method is thus reduced by a
significant margin compared with the similar results in 
Figure~\ref{fig:Sobol_Subsurf}. The efficiency gains of this method
indicate the potential for deployment on problems which would otherwise be
intractable.

\section{Conclusion and future work}\label{sec:Conclusion}

%
%
%
%

We have shown that the feasibility of the standard double-loop approach for GSA
of rare event probabilities can be significantly extended beyond simple
applications. This requires appropriate acceleration methods; in our case,
this is achieved through subset simulation and the choice of a surrogate model
allowing for the analytical calculation of Sobol' indices. This approach is
conceptually simple and does not require the development of new, ad hoc
sensitivity concepts. While we have extended the range of applicability of the
double-loop approach, we acknowledge that more research is needed to deal with
computationally expensive, high-dimensional problems. 

The efficiency of our method crucially depends on working with surrogate models for which sensitivity measures\,---\,here, Sobol' indices\,---\,can be computed cheaply or ``for free"; this clearly and strongly limits the type of GSA  which can be carried out by the approach. More generally, if $q$ is the original QoI and if $\tilde q$ is the resulting QoI for a given surrogate model, more work is needed to understand the relationship between the approximation error $q - \tilde q$ and the resulting GSA error $\mathcal S(q) - \mathcal S(\tilde q)$ where $\mathcal S(\cdot)$ is some sensitivity measure; more explicitly, there may be room for the development ``cheap" surrogate models with moderate approximation errors and \emph{small} GSA errors. 
Additionally, both our sensitivity analysis method as well as surrogate modeling
approach rely on the assumption that the hyper-parameters are independent.
In some cases one might be interested in GSA of rare event probabilities to 
both hyper-parameters and additional parameters in a model that might be 
uncertain and possibly correlated. 
Therefore, another interesting line of inquiry is to consider GSA of rare event 
probability with respect to 
correlated parameters. 
Further study may also include extensions of our approach to other moment-based QoIs (e.g. CDF approximation, skewness, kurtosis) and the use of perturbation-based methods for GSA~\cite{lemaitre2015density} as opposed to considering a discrete set of hyper-parameters.


\paragraph{Acknowledgements}
This research was supported by NSF through grants DMS 1745654 and DMS 1953271 and the US Dept. of Energy (DOE) through Sandia National Laboratories.

	\bibliography{refs}

\begin{thebibliography}{32}
\providecommand{\natexlab}[1]{#1}
\providecommand{\url}[1]{\texttt{#1}}
\expandafter\ifx\csname urlstyle\endcsname\relax
  \providecommand{\doi}[1]{doi: #1}\else
  \providecommand{\doi}{doi: \begingroup \urlstyle{rm}\Url}\fi

\bibitem[spe(2000)]{spe_data}
2001 {SPE} comparative solution project., 2000.
\newblock URL \url{https://www.spe.org/web/csp/datasets/set02.htm}.

\bibitem[Alexanderian(2013)]{Alexanderian13}
Alen Alexanderian.
\newblock On spectral methods for variance based sensitivity analysis.
\newblock \emph{Probability Surveys}, 10:\penalty0 51--68, 2013.

\bibitem[Asmussen and Glynn(2007)]{AsmussenGlynn07}
S{\o}ren Asmussen and Peter~W Glynn.
\newblock \emph{Stochastic simulation: algorithms and analysis}, volume~57.
\newblock Springer Science \& Business Media, 2007.

\bibitem[Au and Beck(2001)]{au2001estimation}
Siu-Kui Au and James~L Beck.
\newblock Estimation of small failure probabilities in high dimensions by
  subset simulation.
\newblock \emph{Probabilistic engineering mechanics}, 16\penalty0 (4):\penalty0
  263--277, 2001.

\bibitem[Beck and Zuev(2017)]{beckzuev}
James~L. Beck and Konstantin~M. Zuev.
\newblock Rare-event simulation.
\newblock In \emph{Handbook of uncertainty quantification. {V}ol. 1, 2, 3},
  pages 1075--1100. Springer, Cham, 2017.

\bibitem[Bhatia and Davis(2000)]{BhatiaDavis00}
Rajendra Bhatia and Chandler Davis.
\newblock A better bound on the variance.
\newblock \emph{The american mathematical monthly}, 107\penalty0 (4):\penalty0
  353--357, 2000.

\bibitem[Blatman and Sudret(2010)]{BlatmanSudret10}
G{\'e}raud Blatman and Bruno Sudret.
\newblock Efficient computation of global sensitivity indices using sparse
  polynomial chaos expansions.
\newblock \emph{Reliab. Eng. Syst. Safe.}, 95\penalty0 (11):\penalty0
  1216--1229, 2010.

\bibitem[Chabridon(2018)]{chabridonthesis}
Vincent Chabridon.
\newblock \emph{Reliability-oriented sensitivity analysis under probabilistic
  model uncertainty--Application to aerospace systems}.
\newblock PhD thesis, Universit{\'e} Clermont Auvergne, 2018.

\bibitem[Chabridon et~al.(2018)Chabridon, Balesdent, Bourinet, Morio, and
  Gayton]{chabridon2018reliability}
Vincent Chabridon, Mathieu Balesdent, Jean-Marc Bourinet, J{\'e}r{\^o}me Morio,
  and Nicolas Gayton.
\newblock Reliability-based sensitivity estimators of rare event probability in
  the presence of distribution parameter uncertainty.
\newblock \emph{Reliability Engineering \& System Safety}, 178:\penalty0
  164--178, 2018.

\bibitem[Cleaves et~al.(2019)Cleaves, Alexanderian, Guy, Smith, and
  Yu]{cleaves2019derivative}
Helen~L Cleaves, Alen Alexanderian, Hayley Guy, Ralph~C Smith, and Meilin Yu.
\newblock Derivative-based global sensitivity analysis for models with
  high-dimensional inputs and functional outputs.
\newblock \emph{SIAM Journal on Scientific Computing}, 41\penalty0
  (6):\penalty0 A3524--A3551, 2019.

\bibitem[Crestaux et~al.(2009)Crestaux, Le~Ma{\i}tre, and
  Martinez]{crestaux2009polynomial}
Thierry Crestaux, Olivier Le~Ma{\i}tre, and Jean-Marc Martinez.
\newblock Polynomial chaos expansion for sensitivity analysis.
\newblock \emph{Reliability Engineering \& System Safety}, 94\penalty0
  (7):\penalty0 1161--1172, 2009.

\bibitem[Dupuis et~al.(2020)Dupuis, Katsoulakis, Pantazis, and
  Rey-Bellet]{dupuis2020sensitivity}
Paul Dupuis, Markos~A Katsoulakis, Yannis Pantazis, and Luc Rey-Bellet.
\newblock Sensitivity analysis for rare events based on r{\'e}nyi divergence.
\newblock \emph{The Annals of Applied Probability}, 30\penalty0 (4):\penalty0
  1507--1533, 2020.

\bibitem[Ehre et~al.(2020)Ehre, Papaioannou, and Straub]{ehre2020framework}
Max Ehre, Iason Papaioannou, and Daniel Straub.
\newblock A framework for global reliability sensitivity analysis in the
  presence of multi-uncertainty.
\newblock \emph{Reliability Engineering \& System Safety}, 195:\penalty0
  106726, 2020.

\bibitem[Fajraoui et~al.(2017)Fajraoui, Marelli, and
  Sudret]{fajraoui2017sequential}
Noura Fajraoui, Stefano Marelli, and Bruno Sudret.
\newblock Sequential design of experiment for sparse polynomial chaos
  expansions.
\newblock \emph{SIAM/ASA Journal on Uncertainty Quantification}, 5\penalty0
  (1):\penalty0 1061--1085, 2017.

\bibitem[Hampton and Doostan(2017)]{hampton2016compressive}
Jerrad Hampton and Alireza Doostan.
\newblock Compressive sampling methods for sparse polynomial chaos expansions.
\newblock In \emph{Handbook of uncertainty quantification}, pages 827--855.
  Springer International Publishing, 2017.

\bibitem[Le~Ma{\^\i}tre and Knio(2010)]{le2010spectral}
Olivier Le~Ma{\^\i}tre and Omar~M Knio.
\newblock \emph{Spectral methods for uncertainty quantification: with
  applications to computational fluid dynamics}.
\newblock Springer Science \& Business Media, 2010.

\bibitem[Lema{\^\i}tre et~al.(2015)Lema{\^\i}tre, Sergienko, Arnaud, Bousquet,
  Gamboa, and Iooss]{lemaitre2015density}
Paul Lema{\^\i}tre, Ekatarina Sergienko, Aur{\'e}lie Arnaud, Nicolas Bousquet,
  Fabrice Gamboa, and Bertrand Iooss.
\newblock Density modification-based reliability sensitivity analysis.
\newblock \emph{Journal of Statistical Computation and Simulation}, 85\penalty0
  (6):\penalty0 1200--1223, 2015.

\bibitem[Li and Xiu(2010)]{li2010evaluation}
Jing Li and Dongbin Xiu.
\newblock Evaluation of failure probability via surrogate models.
\newblock \emph{Journal of Computational Physics}, 229\penalty0 (23):\penalty0
  8966--8980, 2010.

\bibitem[Li et~al.(2011)Li, Li, and Xiu]{li2011efficient}
Jing Li, Jinglai Li, and Dongbin Xiu.
\newblock An efficient surrogate-based method for computing rare failure
  probability.
\newblock \emph{Journal of Computational Physics}, 230\penalty0 (24):\penalty0
  8683--8697, 2011.

\bibitem[Melchers and Beck(2018)]{melchers2018structural}
Robert~E Melchers and Andr{\'e}~T Beck.
\newblock \emph{Structural reliability analysis and prediction}.
\newblock John wiley \& sons, 2018.

\bibitem[Morio(2011)]{morio2011influence}
J{\'e}r{\^o}me Morio.
\newblock Influence of input pdf parameters of a model on a failure probability
  estimation.
\newblock \emph{Simulation Modelling Practice and Theory}, 19\penalty0
  (10):\penalty0 2244--2255, 2011.

\bibitem[Papaioannou et~al.(2015)Papaioannou, Betz, Zwirglmaier, and
  Straub]{papaioannou2015mcmc}
Iason Papaioannou, Wolfgang Betz, Kilian Zwirglmaier, and Daniel Straub.
\newblock {MCMC} algorithms for subset simulation.
\newblock \emph{Probabilistic Engineering Mechanics}, 41:\penalty0 89--103,
  2015.

\bibitem[Peherstorfer et~al.(2017)Peherstorfer, Kramer, and
  Willcox]{peherstorfer2017combining}
Benjamin Peherstorfer, Boris Kramer, and Karen Willcox.
\newblock Combining multiple surrogate models to accelerate failure probability
  estimation with expensive high-fidelity models.
\newblock \emph{Journal of Computational Physics}, 341:\penalty0 61--75, 2017.

\bibitem[Saltelli et~al.(2010)Saltelli, Annoni, Azzini, Campolongo, Ratto, and
  Tarantola]{saltelli2010}
Andrea Saltelli, Paola Annoni, Ivano Azzini, Francesca Campolongo, Marco Ratto,
  and Stefano Tarantola.
\newblock Variance based sensitivity analysis of model output. design and
  estimator for the total sensitivity index.
\newblock \emph{Comput. Phys. Commun.}, 181:\penalty0 259--270, 2010.

\bibitem[Schu{\"e}ller et~al.(2004)Schu{\"e}ller, Pradlwarter, and
  Koutsourelakis]{schueller2004critical}
GI~Schu{\"e}ller, HJ~Pradlwarter, and Phaedon-Stelios Koutsourelakis.
\newblock A critical appraisal of reliability estimation procedures for high
  dimensions.
\newblock \emph{Probabilistic engineering mechanics}, 19\penalty0 (4):\penalty0
  463--474, 2004.

\bibitem[{\v{S}}ehi{\'c} and Karamehmedovi{\'c}(2020)]{vsehic2020estimation}
Kenan {\v{S}}ehi{\'c} and Mirza Karamehmedovi{\'c}.
\newblock Estimation of failure probabilities via local subset approximations.
\newblock \emph{arXiv preprint arXiv:2003.05994}, 2020.

\bibitem[Tong et~al.(2020)Tong, Vanden-Eijnden, and Stadler]{tong2020extreme}
Shanyin Tong, Eric Vanden-Eijnden, and Georg Stadler.
\newblock Extreme event probability estimation using pde-constrained
  optimization and large deviation theory, with application to tsunamis.
\newblock \emph{arXiv preprint arXiv:2007.13930}, 2020.

\bibitem[Ullmann and Papaioannou(2015)]{ullmann2015multilevel}
Elisabeth Ullmann and Iason Papaioannou.
\newblock Multilevel estimation of rare events.
\newblock \emph{SIAM/ASA Journal on Uncertainty Quantification}, 3\penalty0
  (1):\penalty0 922--953, 2015.

\bibitem[van~den Berg and Friedlander(2019)]{spgl1site}
E.~van~den Berg and M.~P. Friedlander.
\newblock {SPGL1}: A solver for large-scale sparse reconstruction, December
  2019.
\newblock https://friedlander.io/spgl1.

\bibitem[Wang et~al.(2021)Wang, Li, Liu, and Zhou]{wang2021global}
Pan Wang, Chunyu Li, Fuchao Liu, and Hanyuan Zhou.
\newblock Global sensitivity analysis of failure probability of structures with
  uncertainties of random variable and their distribution parameters.
\newblock \emph{Engineering with Computers}, pages 1--19, 2021.

\bibitem[Wang and Jia(2020)]{wang2020augmented}
Zhenqiang Wang and Gaofeng Jia.
\newblock Augmented sample-based approach for efficient evaluation of risk
  sensitivity with respect to epistemic uncertainty in distribution parameters.
\newblock \emph{Reliability Engineering \& System Safety}, 197:\penalty0
  106783, 2020.

\bibitem[Zuev et~al.(2012)Zuev, Beck, Au, and Katafygiotis]{zuev2012bayesian}
Konstantin~M Zuev, James~L Beck, Siu-Kui Au, and Lambros~S Katafygiotis.
\newblock Bayesian post-processor and other enhancements of subset simulation
  for estimating failure probabilities in high dimensions.
\newblock \emph{Computers \& structures}, 92:\penalty0 283--296, 2012.

\end{thebibliography}

\end{document}